\documentclass{amsproc}

\def\g{\gamma}

\def\ep{\epsilon}
\def\a{\alpha}
\def\b{\beta}
\def\d{\delta}
\def\be{\begin{equation}}
\def\ee{\end{equation}}
\def\bea{\begin{eqnarray}}
\def\eea{\end{eqnarray}}
\def\ba{\begin{array}}
\def\ea{\end{array}}

\def\bx{{\bf x}}
\def\by{{\bf y}}
\def\bl{{\bf l}}
\def\bc{{\bf c}}

\def\bff{{\bf f}}
\def\bg{{\bf g}}

\def\ra{\rightarrow}
\def\R{{\mathbb R}}

\def\C{{\mathbb C}}
\def\g{\mathcal{G}}

\newtheorem{theorem}{Theorem}[section]
\newtheorem{lemma}[theorem]{Lemma}

\theoremstyle{definition}

\theoremstyle{remark}

\numberwithin{equation}{section}

\newcommand{\abs}[1]{\lvert#1\rvert}

\newfont{\Bb}{msbm8 scaled\magstep{1}}
\newcommand{\rc}{\mbox{\Bb R}}

\usepackage{graphics}

\begin{document}

\title[Optimization methods in direct and inverse scattering]
{Optimization methods in direct and inverse scattering
}

\author{Alexander G. RAMM}
\address{ Department of Mathematics\\ Kansas State University\\
Manhattan, Kansas 66506-2602, USA}
\email{ramm@math.ksu.edu}

\author{Semion GUTMAN}
\address{Department of Mathematics\\ University of Oklahoma\\ Norman,
 OK 73019, USA}
\email{sgutman@ou.edu}

\subjclass{Primary 78A46, 65N21, Secondary 35R30}

\begin{abstract}
In many Direct and Inverse Scattering problems
 one has to use a parameter-fitting procedure,
because analytical inversion procedures are often not available. In this
paper a variety of such methods is presented with a discussion
of theoretical and computational issues.

The problem of finding  small subsurface inclusions
from surface scattering data is stated and investigated. This Inverse
Scattering problem is reduced
to  an optimization problem, and solved by
the Hybrid
Stochastic-Deterministic minimization algorithm.
A similar approach is used to determine layers in a
particle from the scattering data.

The Inverse potential scattering problem is described and its solution
based on a parameter fitting procedure is presented for the case
of spherically symmetric potentials and fixed-energy phase shifts as the
scattering data.
 The central feature of the minimization algorithm here is
the Stability Index Method. This general approach estimates the size of
the minimizing sets, and gives a practically useful stopping criterion for
global minimization algorithms.

The  3D inverse scattering problem with
fixed-energy data is discussed. Its solution by the Ramm's method is
described. The cases of exact and noisy discrete data
are considered. Error estimates for the inversion algorithm are given in
both cases of exact and noisy data. Comparison of the Ramm's inversion
method with the inversion based on the Dirichlet-to-Neumann map is given
and it is shown that there are many more numerical difficulties in the
latter method than in the Ramm's method.

An Obstacle Direct Scattering problem is
treated by a novel
Modified Rayleigh Conjecture (MRC) method.  MRC's performance is compared
favorably to the well known Boundary Integral Equation Method, based on
the properties of the single and double-layer potentials. A special
minimization procedure allows one to inexpensively compute 
scattered
fields for 2D and 3D obstacles having smooth as well as
nonsmooth surfaces.

A new Support Function Method (SFM) is used for Inverse Obstacle
Scattering problems. The SFM can work with limited
data. It can also be used for Inverse scattering problems with {\it
unknown scattering conditions on its boundary (e.g. soft, or hard
scattering)}. Another method for Inverse scattering problems, the
Linear Sampling Method (LSM), is analyzed.
Theoretical and computational difficulties in using this method are
pointed out.

\end{abstract}

\maketitle

Key words and phrases: inverse and direct scattering,
 optimization, Modified Rayleigh Conjecture, support function method,
stability index, Ramm's method, small inhomogeneities, linear sampling
method.

\tableofcontents

\section{Introduction.}

Suppose that an acoustic or electromagnetic wave encounters an
inhomogeneity and, as a consequence, gets scattered. The
problem of finding the scattered wave assuming the knowledge of the
inhomogeneity (penetrable or not) is the Direct Scattering
problem. An impenetrable inhomogeneity is also called an obstacle.
On the other hand, if the scattered wave is known at some points
outside an inhomogeneity,  then we are
faced with the Inverse Scattering problem, the goal of which is
to identify this inhomogeneity, see
\cite{r190,rammb2,r313,coltonkress,cc}

Among a variety of methods available to handle such problems few
provide a mathematically justified algorithm. In many cases one
has to use a parameter-fitting procedure, especially for inverse
scattering problems, because the analytical inversion procedures
are often not available. An important part of such a procedure is
an efficient global optimization method, see
\cite{floudas1,floudas2, horst1, horst2, pardalos, ru}.

 The general scheme for parameter-fitting procedures is simple: one has a
relation
$B(q)=A$, where $B$ is some operator, $q$ is an unknown function,
and $A$ is the data. In inverse scattering problems $q$ is an
unknown potential,
and $A$ is the known scattering amplitude.
 If $q$ is sought in a finite-parametric family of functions, then
$q=q(x, p)$, where $p=(p_1,....,p_n)$ is a parameter. The parameter
is found by solving a global minimization problem:
$\Phi[B(q(x,p))-A]=\min$, where $\Phi$ is some positive functional,
and $q\in Q$, where $Q$ is an admissible set of $q$.
   In practice the above problem often has many local minimizers,
and the global minimizer is not necessarily unique.
In \cite{rammb2} and \cite{r317} some functionals $\Phi$ are constructed
which have unique global minimizer, namely, the solution to inverse
scattering problem, and the global minimum is zero.

 Moreover, as a rule, the data $A$ is known with some error. Thus
$A_\delta$ is known, such that $||A-A_\delta||<\delta$. There are
no stability estimates which would show how the global minimizer
$q(x,p_{opt})$
is perturbed when the data $A$ are replaced by the perturbed data
$A_\delta$. In fact, one can easily construct examples showing that
there is no stability of the global minimizer with respect to small
errors in the data, in general.

 For these reasons  there is no guarantee that the
parameter-fitting procedures would yield a solution to the inverse problem
with a guaranteed accuracy.
  However, overwhelming majority of practitioners are using
parameter-fitting procedures. In dozens of published papers
the results obtained by various parameter-fitting procedures look quite
good. The explanation, in most of the cases is simple: the authors know
the answer beforehand, and it is usually not difficult to
parametrize the unknown function so that the exact solution is well
approximated by a function from a finite-parametric family, and since
the authors know a priori the exact answer, they may choose
numerically the values of the parameters which yield a good approximation
of the exact solution.
 {\it When can one rely on the results obtained by
parameter-fitting procedures? Unfortunately, there is no rigorous and
complete answer to this question, but some recommendations are given in
Section 4. }

In this paper the authors present their recent results which are based on
specially designed parameter-fitting procedures. Before describing them,
let us mention that usually in a numerical solution of an inverse
scattering problem one uses a regularization procedure, e.g. a variational
regularization,  spectral cut-off, iterative regularization,
DSM (the dynamical systems method), quasi-solutions, etc, see e.g.
\cite{r470}, \cite{r456}.
This general theoretical framework is well
established in the theory of ill-posed problems, of which the inverse
scattering problems represent an important class. This framework is needed
to achieve a stable method for assigning a solution to an ill-posed
problem, usually set in an infinite dimensional space. The goal of
this paper is to present optimization algorithms already in a finite
dimensional setting of a Direct or Inverse scattering problem.

In Section 2 the problem of finding  small subsurface inclusions
from
surface scattering data is investigated (\cite{R4}). This (geophysical)
Inverse
Scattering problem is reduced
to  an optimization problem. This problem is solved by
the Hybrid
Stochastic-Deterministic minimization algorithm (\cite{r399}). It is based
on a genetic
minimization
algorithm ideas for its random (stochastic) part, and a deterministic
minimization without derivatives used for the local minimization part.

In Section 3 a similar approach is used to determine layers in a
particle subjected to acoustic or electromagnetic waves. The global
minimization algorithm uses Rinnooy Kan and Timmer's
Multilevel Single-Linkage Method for its stochastic part.

In Section 4 we discuss an Inverse potential scattering problem
appearing in a quantum mechanical description of particle scattering
experiments. The central feature of the minimization algorithm here is
the Stability Index Method (\cite{r435}). This general approach estimates
the size of
the minimizing sets, and gives a practically useful stopping criterion for
global
minimization algorithms.

In Section 5 Ramm's method for solving 3D inverse scattering problem with
fixed-energy data is presented following \cite{r474}, see also
\cite{r425}. The cases
of exact and noisy discrete data
are considered. Error estimates for the inversion algorithm are given in
both cases of exact and noisy data. Comparison of the Ramm's inversion
method with the inversion based on the Dirichlet-to-Neumann map is given
and it is shown that there are many more numerical difficulties in the
latter method than in Ramm's method.

 In Section 6 an Obstacle Direct Scattering problem is
treated by a novel
Modified Rayleigh Conjecture (MRC) method. It was introduced in \cite{r430}
and applied in \cite{r437,r475} and \cite{r461}. MRC's performance is compared
favorably to the well known Boundary Integral Equation Method, based on
the properties of the single and double-layer potentials. A special
minimization procedure allows us to inexpensivly compute scattered
fields for several 2D and 3D obstacles having smooth as well as
nonsmooth surfaces.

In Section 7 a new Support Function Method (SFM) is used to determine
the location of an obstacle (cf \cite{r56}, \cite{r190}, \cite{r453}).
Unlike other methods, the SFM can work with
limited data. It can also be used for Inverse scattering problems
with {\it unknown scattering conditions on its boundary (e.g. soft, or
hard scattering)}.

Finally, in Section 8, we present an analysis of another popular method
for Inverse scattering problems, the Linear Sampling Method (LSM), and
show that both theoretically and computationally the method fails in
many aspects.

\section{Identification of small subsurface inclusions.}

\subsection{Problem description.}
In many applications it is desirable to find small
inhomogeneities from surface scattering data. For example, such a problem
arises in ultrasound mammography, where small inhomogeneities
are cancer cells.
Other examples include the problem of finding small holes and cracks
in metals and other materials, or the mine detection.
The scattering theory for small scatterers originated in the
 classical works of Lord Rayleigh (1871). Rayleigh understood that the
basic contribution to the scattered field in the far-field zone comes
from the dipole radiation, but did not give methods for calculating this
radiation. Analytical formulas for calculating the polarizability
tensors for homogeneous bodies of arbitrary shapes were derived in
\cite{r190} (see also references therein). These formulas allow one
to calculate the $S$-matrix for scattering of acoustic and
electromagnetic waves by small bodies of arbitrary shapes with arbitrary
accuracy.
  Inverse scattering problems for small bodies are considered in \cite{r1}
and \cite{r313}.
In \cite{R4} the problem of identification of small subsurface inhomogeneities
from surface data was posed and its possible applications were discussed.

In the context of a geophysical problem, let  $y\in\R^3$ be
a point source of monochromatic acoustic waves on the
surface of the earth. Let $u(x,y,k)$ be the acoustic pressure at
a point $x\in\R^3$, and $k>0$ be the wavenumber. The governing equation  for the acoustic wave propagation is:

\begin{equation}\label{subs_e1}
\left[\nabla^2+k^2+k^2v(x)\right] u=-\delta(x-y)\ \text{in}\ \R^3,
\end{equation}
where $x=(x_1,x_2,x_3),\; v(x)$
is the inhomogeneity in the velocity profile, and $u(x,y,k)$ satisfies the radiation condition at infinity,
i.e. it decays sufficiently fast as $|x|\ra\infty$.

Let us assume that $v(x)$ is a bounded function vanishing outside of the domain
$D=\cup_{m=1}^M D_m $
which is the union of $M$ small nonintersecting domains $D_m$,
all of them are located in the lower half-space
$\R^3_{-}=\{x:x_3<0\}$. Smallness is understood in the sense
$k\rho\ll 1$, where $\rho:=\frac 12 \max_{1\leq m\leq
M}\{\text{diam}\,D_m\}$,
and diam $D$ is the diameter of the domain $D$. Practically
$k\rho\ll 1$ means that $k\rho<0.1$. In some cases $k\rho<0.2$ is
sufficient
for obtaining acceptable numerical results. The background velocity
in (\ref{subs_e1}) equals to 1, but we can consider the case of fairly
general background velocity \cite{r313}.

Denote $\tilde z_m$ and $\tilde v_m$ the position of the center of gravity
 of $D_m$, and the total intensity of the $m\text{-th}$
inhomogeneity $\tilde v_m:=\int_{D_m}v(x)dx$. Assume that
$\tilde v_m\not= 0$. Let $P$ be the equation of the surface of the earth:
\begin{equation}
P:=\{x=(x_1,x_2,x_3)\in\R^3 : x_3=0\}.
\end{equation}

The inverse problem to be solved is:

{\bf IP:} {\it Given $u(x,y,k)$ for all source-detector pairs $(x,y)$ on
$P$ at a fixed $k>0$, find
the number $M$ of small inhomogeneities,
the positions $\tilde z_m$ of the inhomogeneities,
and their intensities $\tilde v_m$.
}

Practically, one assumes that
 a fixed wavenumber $k>0$, and $J$ source-detector pairs $(x_j,y_j),
j=1,2,...,J,$ on $P$ are known together with the acoustic pressure
measurements $u(x_j, y_j,k)$. Let \begin{equation}
g(x,y,k):=\frac{\exp(ik|x-y|)}{4\pi|x-y|},\quad x,y\in P, \end{equation}
\begin{equation} G_j(z):=G(x_j,y_j,z):=g(x_j,z,k)g(y_j,z,k),\quad
x_j,y_j\in P,\; z\in \R^3_-, \end{equation} \begin{equation}
f_j:=\frac{u(x_j,y_j,k)-g(x_j,y_j,k)}{k^2}, \end{equation} and
\begin{equation}\label{subs_e6}
 \Phi(z_1,\dots,z_M,\, v_1,\dots,v_M):=
   \sum^J_{j=1}\left| f_j-\sum^M_{m=1} G_j(z_m)v_m\right|^2.
\end{equation}

The proposed method for solving the (IP) consists of finding the
global minimizer of function (\ref{subs_e6}). This minimizer
$(\tilde z_1,\dots,\tilde z_M,\, \tilde v_1,\dots,\tilde v_M)$
gives the estimates  of the
positions $\tilde z_m$ of the small inhomogeneities and their intensities
$\tilde v_m$. See \cite{R4} for a justification of this approach.

The function $\Phi$ depends on $M$ unknown points
$z_m\in{\R}^3_-$, and $M$ unknown parameters $v_m$,
$1\leq m\leq M$. The    number $M$ of the small
inhomogeneities is also unknown, and its determination is a part of the
minimization problem.

\subsection{Hybrid Stochastic-Deterministic Method}

Let the inhomogeneities be located within the box
\begin{equation}
B=\{(x_1,x_2,x_3)\ :-a<x_1<a,\ -b<x_2<b,\ 0<x_3<c\}\,,
\end{equation}
and their intensities satisfy
\begin{equation}\label{surf_e8}
0\leq v_m\leq v_{max}\,.
\end{equation}
The box is located above the earth surface for a computational
convenience.

Then, given the location of the points $z_1, z_2, \dots, z_M$, the
minimum of $\Phi$  with respect to the intensities $v_1, v_2,
\dots, v_M$ can be found by minimizing the resulting quadratic function in (\ref{subs_e6}) over the region
satisfying (\ref{surf_e8}). This can be done using normal equations for (\ref{subs_e6}) and
projecting the resulting point back onto the region defined by (\ref{surf_e8}).
Denote the result of this minimization by $\tilde \Phi$, that is

\begin{equation}\label{surf_e9}
\begin{split}
\tilde \Phi(z_1, z_2, \dots, z_M)& =\min\{\Phi(z_1, z_2, \dots, z_M,
v_1, v_2, \dots, v_M)\ :\\
&\ 0\leq v_m\leq v_{max}\,,\quad 1\leq m \leq M\}
\end{split}
\end{equation}

Now the original minimization problem for $\Phi(z_1, z_2, \dots, z_M,
v_1, v_2, \dots, v_M) $ is reduced to the $3M$-dimensional
constrained minimization for $\tilde \Phi(z_1, z_2, \dots, z_M) $:
\begin{equation}\label{surf_e10}
\tilde \Phi(z_1, z_2, \dots, z_M)=\min \,,\quad z_m\in B\,,\quad 1\leq m
\leq M\,.
\end{equation}

Note, that the dependency of $\tilde\Phi$ on its $3M$ variables (the
coordinates of the points $z_m$) is highly nonlinear. In particular,
this dependency is complicated by the computation of the minimum in
(\ref{surf_e9}) and the consequent projection onto the admissible set $B$. Thus,
an analytical computation of the gradient of $\tilde\Phi$ is not
computationally efficient. Accordingly, the Powell's quadratic minimization
method was used to find local minima. This method uses a special procedure to numerically approximate the
gradient, and it can be shown to exhibit the same type of quadratic
convergence as conjugate gradient type methods (see \cite{bre}).

\begin{figure}[tb]
\vspace{5pc}
\includegraphics*{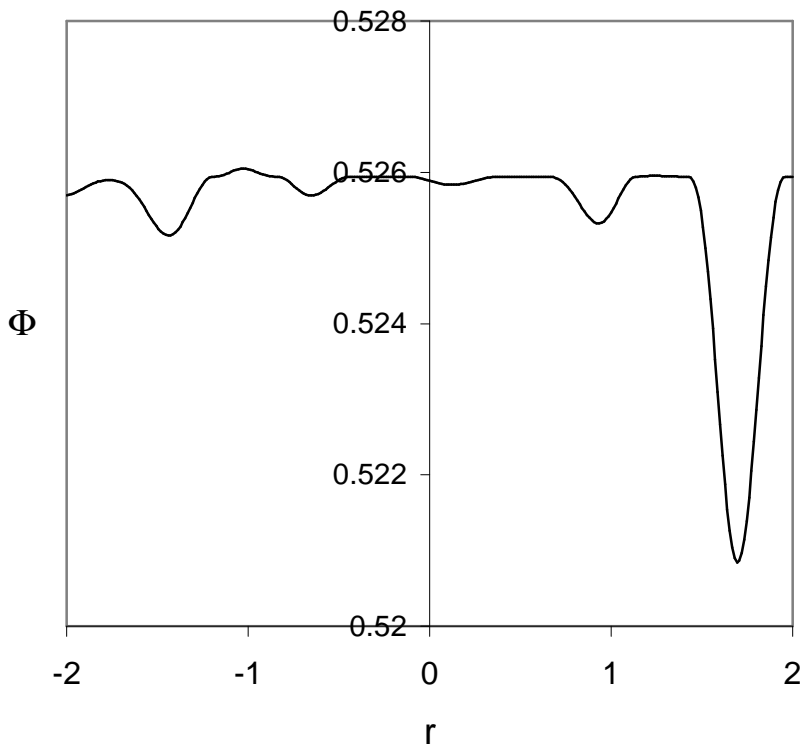}
\caption{Objective
function $\tilde\Phi(z_r, z_2,\tilde z_3,\tilde z_4,\tilde z_5,\tilde
z_6)\,,\quad -2\leq r\leq 2$ }\label{subs_fig1}
\end{figure}

In addition,
{\it the exact number of the original inhomogeneities $M_{orig}$ is
unknown,
and its estimate is a part of the inverse problem.} In the HSD algorithm
described below this task is accomplished by taking the initial number
$M$ sufficiently large, so that
\begin{equation}
M_{orig}\leq M\,,
\end{equation}
which, presumably, can be estimated from physical considerations. After
all, our goal  is to find only the strongest inclusions, since
the weak ones cannot be distinguished from  background noise. The
Reduction Procedure (see below) allows the algorithm to seek the minimum
of $\tilde\Phi$ in a lower dimensional subsets of the admissible set
$B$, thus finding the estimated number of inclusions $M$. Still another
difficulty in the minimization is a large number of local minima of
$\tilde \Phi$. This phenomenon is well known for  objective
functions arising in various inverse problems, and we illustrate this
point
in Figure \ref{subs_fig1}.

For example, let $M_{orig}=6$, and the coordinates of the inclusions, and their
intensities
$(\tilde z_1,  \dots, \tilde z_6,
\tilde v_1,  \dots, \tilde v_6)$
be as in Table \ref{subs_table1}. Figure \ref{subs_fig1} shows the values of the
function $\tilde\Phi(z_r, z_2,\tilde z_3,\tilde z_4,\tilde z_5,\tilde
z_6)$, where
\[
z_r=(r,0,0.520)\,,\quad -2\leq r\leq 2
\]
and
\[z_2=(-1,0.3,0.580)\,.
\]
The plot shows multiple local minima and almost flat regions.

A direct application of a gradient type method to such a function would
result in finding a local minimum, which may or may not be the sought
global one. In the example above, such a method would usually be trapped
in a local minimum located at $r=-2,\ r=-1.4,\ r=-0.6,\ r=0.2$ or $r=0.9$,
and the desired global minimum at $r=1.6$ would be found only for a
sufficiently close initial guess $1.4<r<1.9$. Various global minimization
methods are known (see below), but we found that an efficient way to
accomplish the minimization task for this Inverse Problem was to design a
new method (HSD) combining both the stochastic and the deterministic
approach to the global minimization. Deterministic minimization algorithms
with or without the gradient computation, such as the conjugate gradient
methods, are known to be efficient (see \cite{bre,des,pol,jac}), and
\cite{ru}.
However, the initial guess should be chosen sufficiently close to the
sought minimum.  Also such algorithms tend to be
trapped at a local minimum, which is not necessarily close to a global
one. A new deterministic method is proposed in \cite{bar} and \cite{bpr},
which is quite efficient according to \cite{bpr}. On the other hand,
various stochastic minimization algorithms, e.g. the simulated annealing
method \cite{kgv,kir}, are more likely to find a global minimum, but their
convergence can be very slow. We have tried a variety of minimization
algorithms to find an acceptable minimum of $\tilde \Phi$. Among them were
the Levenberg-Marquardt Method, Conjugate Gradients, Downhill Simplex, and
Simulated Annealing Method. None of them produced consistent satisfactory
results.

Among minimization methods combining random and deterministic searches we
mention Deep's method \cite{Deep} and a variety of clustering methods
\cite{rkt1}, \cite{rkt2}. An application of these methods to the particle
identification using light scattering is described in \cite{zub}. The
clustering methods are quite robust (that is, they consistently find
global extrema) but, usually, require a significant computational effort.
One such method is described in the next section on the identification of
layers in a multilayer particle.  The HSD method is a combination of a
reduced sample random search method with certain ideas from Genetic
Algorithms (see e.g. \cite{hh}). It is very efficient and seems especially
well suited for low dimensional global minimization. Further research is
envisioned to study its properties in more detail, and its applicability
to other problems.

 The steps of the Hybrid Stochastic-Deterministic
 (HSD) method are outlined below.
 Let us call a collection of
$M$ points ( inclusion's centers)  $\{z_1,z_2,...,z_M\},\; z_i\in B$ a {\it configuration} $Z$.
Then the minimization problem (\ref{surf_e10}) is the minimization of
the objective function $\tilde\Phi$ over the set of all
configurations.

For clarity, let $P_0=1,\; \epsilon_s=0.5,\; \epsilon_i=0.25,\; \ep_d=0.1$,
be the same values as the ones used in numerical computations in the next section.

Generate a random configuration $Z$. Compute the best fit intensities
$v_i$ corresponding to this configuration. If $v_i>v_{max}$, then let
$v_i:=v_{max}$. If $v_i<0$, then let $v_i:=0$.
 If $\tilde\Phi(Z)<P_0\epsilon_s$,
then this
configuration is a preliminary candidate for the initial
guess of a deterministic minimization method (Step 1).

Drop
the points $z_i\in Z$ such that $v_i<v_{max}\epsilon_i$. That is,
the inclusions with small intensities are eliminated (Step 2).

If two points $z_k, z_j\in Z$ are too close to each other, then replace
them with one point of a combined intensity (Step 3).

After completing steps 2 and 3 we would be left with $N\leq M$ points
$z_1,z_2,...,z_N$ (after a re-indexing) of the original configuration $Z$. Use this reduced
configuration $Z_{red}$ as the starting point for the deterministic
restraint minimization in the $3N$ dimensional space (Step 4).  Let the
resulting minimizer be $\tilde Z_{red}=(\tilde z_1,...,\tilde z_N)$. If the
value of the objective function $\tilde \Phi(\tilde Z_{red})<\ep$, then we are
done: $\tilde Z_{red}$ is the sought configuration containing $N$ inclusions.
If $\tilde \Phi(\tilde Z_{red})\geq\ep$, then the iterations should continue.

To continue the iteration, randomly generate $M-N$ points in $B$ (Step 5). Add them
to the reduced configuration $\tilde Z_{red}$. Now we have a new full configuration $Z$, and
the iteration process can continue (Step 1).

This entire iterative process is repeated $n_{max}$ times, and the best
configuration is declared to represent the sought inclusions.

\subsection{Hybrid Stochastic-Deterministic (HSD) Method.}
Let $P_0$, $T_{max}$, $n_{max}$, $ \epsilon_s,\; \epsilon_i,\; \epsilon_d$, and
$\epsilon$ be
positive numbers. Let a positive integer $M$ be larger than the expected
number of inclusions.
 Let $N=0$.

\begin{enumerate}

\item  Randomly generate
$M-N$ additional points $z_{N+1},\dots,z_M\in B$ to obtain a full
configuration $Z=(z_1,\dots,z_M)$. Find the best fit intensities $v_i,\;
i=1,2,...,M$. If $v_i>v_{max}$, then let
$v_i:=v_{max}$. If $v_i<0$, then let $v_i:=0$.
Compute $P_s=\tilde\Phi(z_1,z_2 \dots, z_M)$.
If $P_s<P_0\epsilon_s$ then go to step 2, otherwise repeat
step 1.

\item  Drop all the points with the intensities $v_i$
satisfying $v_i<v_{max}\epsilon_i$.
Now only $N\leq M$ points $z_1,z_2 \dots, z_N$ (re-indexed) remain in
the configuration $Z$.

\item If any two points $z_m, z_n$ in the above configuration satisfy
$|z_m-z_n|<\epsilon_dD$, where $D=diam(B)$, then eliminate point $z_n$,
change the intensity of point $z_m$ to $v_m+v_n$, and assign $N:=N-1$.
This step is repeated until no further reduction in $N$ is possible.
Call the resulting reduced configuration with $N$ points by $Z_{red}$.

\item Run a constrained deterministic minimization of $\tilde\Phi$
in $3N$ variables, with the
initial guess $Z_{red}$.
Let the minimizer be
$\tilde Z_{red}=(\tilde z_1,\dots,\tilde z_N)$.
If $P=\tilde\Phi(\tilde z_1,\dots,\tilde z_N)$ $<\epsilon$, then save this
configuration,
and go to step 6, otherwise let $P_0=P$, and
proceed to the next step 5.

\item Keep intact N points $\tilde z_1,\dots,\tilde z_N$.
If the number of random
configurations has exceeded $T_{max}$ (the maximum
number of random tries), then save the configuration
 and go to step 6, otherwise go to step 1, and use these $N$ points there.

\item  Repeat steps 1 through 5 $n_{max}$ times.

\item Find the configuration among the above $n_{max}$ ones, which gives the smallest
value to $\tilde\Phi$. This is the best fit.

\end{enumerate}

The Powell's minimization method
(see \cite{bre} for a detailed description) was used for the deterministic part, since
this method does not need
gradient computations, and it converges quadratically near quadratically
shaped minima. Also, in step 1, an idea from the Genetic Algorithm's approach \cite{hh} is
implemented by keeping only the strongest representatives of the
population, and allowing a mutation for the rest.

\subsection{Numerical results}

The algorithm was tested on a variety
of configurations. Here
we present the results of just
two typical numerical experiments illustrating the performance of the
method.
In both experiments the box $B$ is taken to be
$$B=\{(x_1,x_2,x_3)\ :-a<x_1<a,\ -b<x_2<b,\ 0<x_3<c\}\,,$$
with $a=2,\ b=1,\ c=1$. The wavenumber $k=5$, and the effective
intensities $v_m$ are in the range from $0$ to $2$.
The values of the parameters  were chosen as follows
$$P_0=1\,,T_{max}=1000,\,\epsilon_s=0.5,\,\epsilon_i=0.25,\,\epsilon_d=0.1\,,
\epsilon=10^{-5}\,,
n_{max}=6$$
In both cases we searched for the same 6 inhomogeneities with the
 coordinates $x_1,x_2,x_3$ and the intensities $v$ shown in Table \ref{subs_table1}.

\begin{table}
\caption{Actual inclusions.}
\begin{tabular}{|r|r|r|r|r|}
\hline
Inclusions &    $\qquad   x_1$ & $\qquad x_2$ & $\qquad x_3$ & $\qquad v$ \\
\hline
1 & 1.640 & -0.510 & 0.520 & 1.200 \\
2 & -1.430 & -0.500 & 0.580 & 0.500 \\
3 & 1.220 & 0.570 & 0.370 & 0.700 \\
4 & 1.410 & 0.230 & 0.740 & 0.610 \\
5 & -0.220 & 0.470 & 0.270 & 0.700 \\
6 & -1.410 & 0.230 & 0.174 & 0.600 \\
\hline
\end{tabular}\label{subs_table1}

\end{table}

Parameter $M$ was set to 16, thus the only information on the number of
inhomogeneities given to the algorithm was that their number does not
exceed 16. This number was chosen to keep the computational time within
reasonable limits. Still another consideration for the number $M$ is the
aim of the algorithm to find the presence of the most influential
inclusions, rather then all inclusions, which is usually impossible in
the presence of noise and with the limited amount of data.

\begin{figure}[tb]
\vspace{5pc}
\includegraphics*{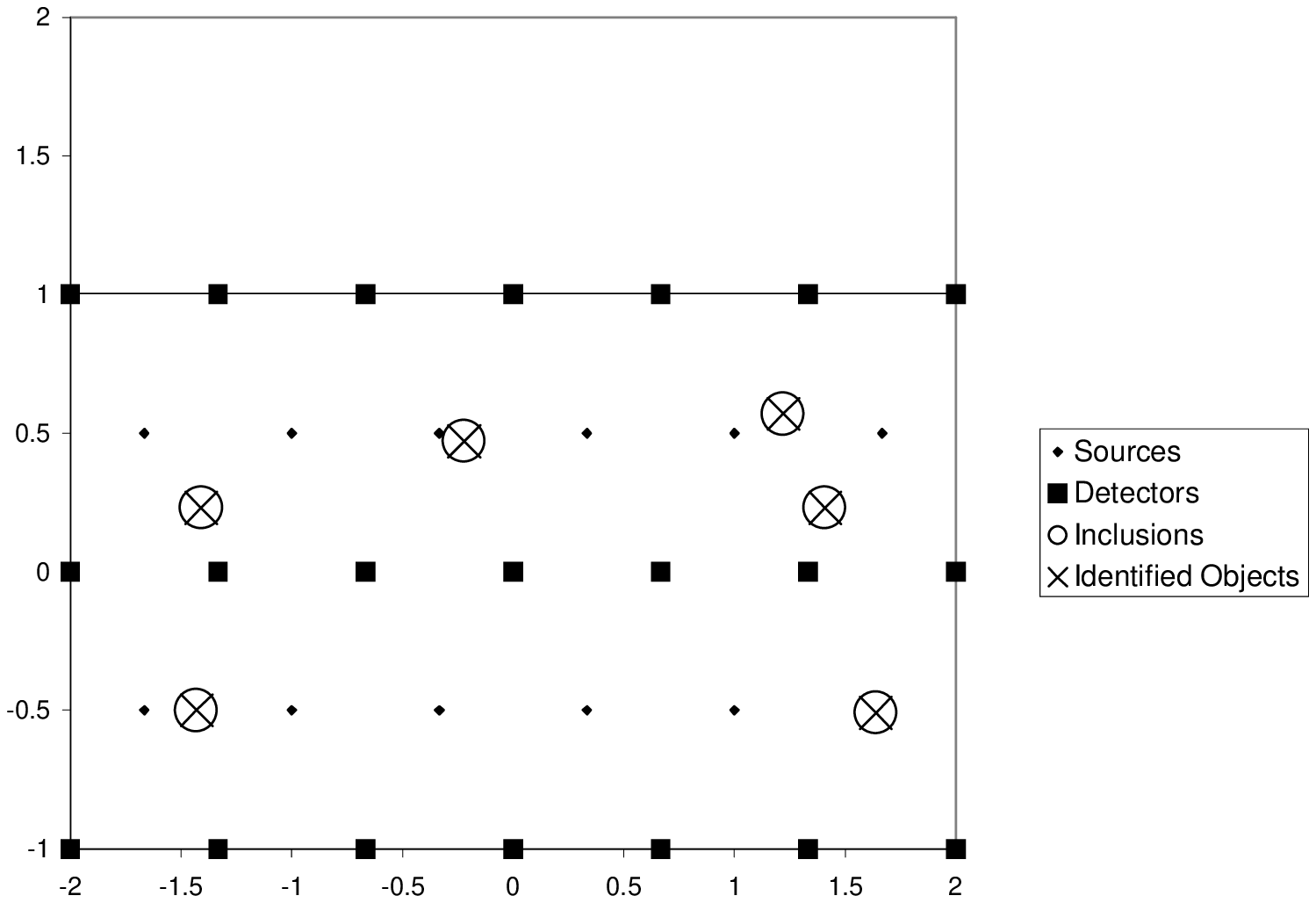}
\caption{Inclusions and Identified objects for subsurface particle identification,
Experiment 1, $\delta=0.00$. $x_3$ coordinate is not
shown.}\label{subs_fig2}
\end{figure}

{\bf Experiment 1.}
In this case we used 12 sources and 21 detectors, all on the surface
$x_3=0$.
The sources were positioned at $\{(-1.667+0.667i,-0.5+1.0j,0),\
i=0,1,\dots,5,\,j=0,1\}$, that is 6 each along two lines $x_2=-0.5$ and
$x_2=0.5$.
The detectors were positioned at $\{(-2+0.667i,-1.0+1.0j,0),\
i=0,1,\dots,6,\,j=0,1,2\}$, that is seven detectors along each of the three lines
$x_2=-1,\, x_2=0$ and
$x_2=1$.
This corresponds to a mammography
search, where the detectors and the sources are placed above the search
area. The results for noise level $\d=0.00$ are shown in 
Figure \ref{subs_fig2} and Table \ref{subs_table2}. The results for noise 
level $\d=0.05$ are shown in Table \ref{subs_table3}.

\begin{table}
\caption{Experiment 1. Identified inclusions, no noise, $\delta=0.00$.}
\begin{tabular}{|r|r|r|r|}
\hline
 $\qquad   x_1$ & $\qquad x_2$ & $\qquad x_3$ & $\qquad v$ \\
\hline
 1.640 & -0.510 & 0.520 & 1.20000 \\
 -1.430 & -0.500 & 0.580 & 0.50000 \\
  1.220 & 0.570 & 0.370 & 0.70000 \\
 1.410 & 0.230 & 0.740 & 0.61000 \\
 -0.220 & 0.470 & 0.270 & 0.70000 \\
-1.410 & 0.230 & 0.174 & 0.60000 \\
\hline
\end{tabular}\label{subs_table2}

\end{table} 

\begin{table}
\caption{Experiment 1. Identified inclusions, $\delta=0.05$.}
\begin{tabular}{|r|r|r|r|}
\hline
 $\qquad   x_1$ & $\qquad x_2$ & $\qquad x_3$ & $\qquad v$ \\
\hline
 1.645 & -0.507 & 0.525 & 1.24243 \\
 1.215 & 0.609 & 0.376 & 0.67626 \\
  -0.216 & 0.465 & 0.275 & 0.69180 \\
-1.395 & 0.248 & 0.177 & 0.60747 \\
\hline
\end{tabular}\label{subs_table3}

\end{table} 

\begin{figure}[tb]
\vspace{5pc}
\includegraphics*{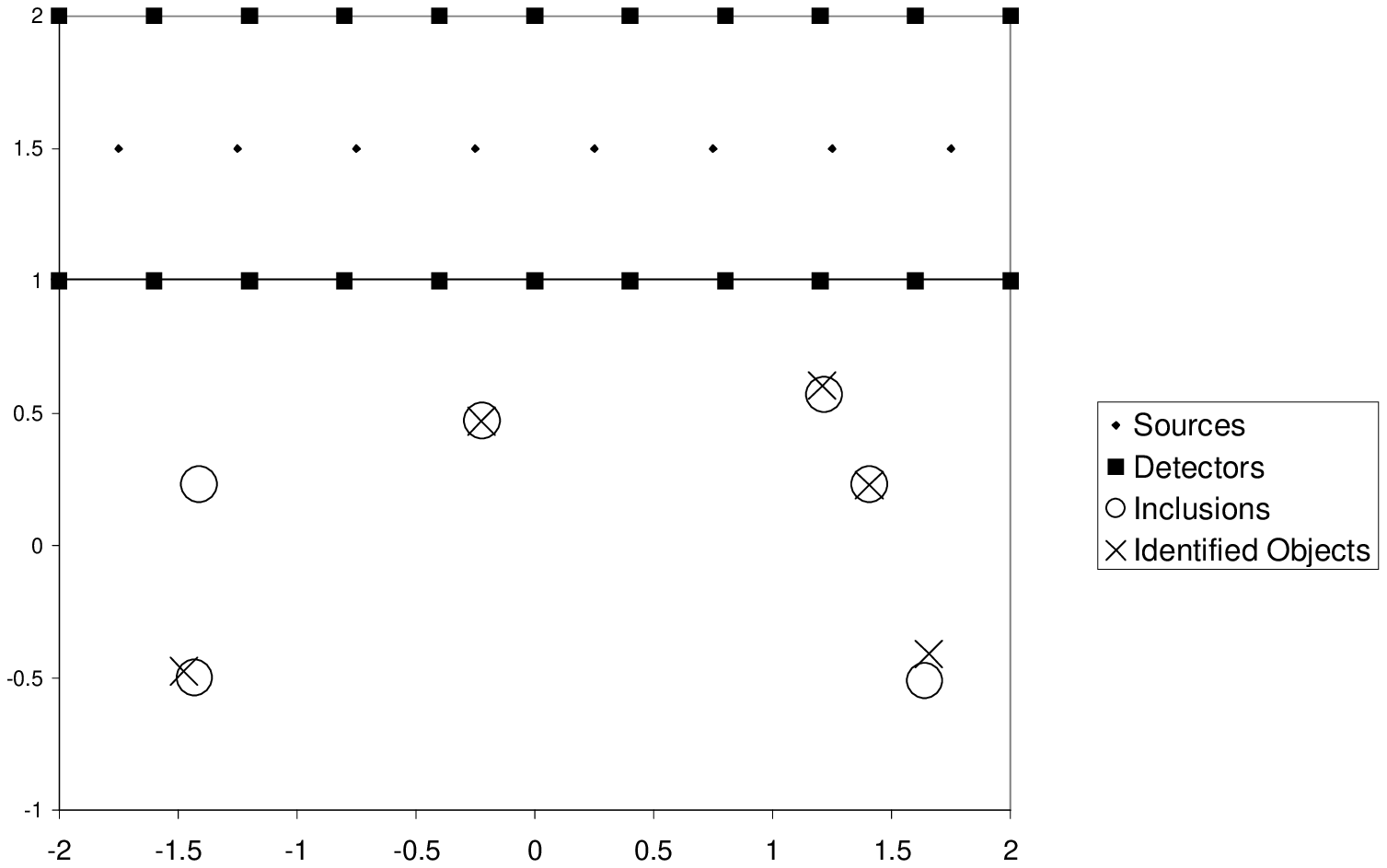}
\caption{Inclusions and Identified objects for for subsurface particle
identification,
Experiment 2, $\delta=0.00$. $x_3$ coordinate is not
shown.}\label{subs_fig3}
\end{figure}

{\bf Experiment 2.}
In this case we used 8 sources and 22 detectors, all on the surface
$x_3=0$.
The sources were positioned at $\{(-1.75+0.5i,1.5,0),\
i=0,1,\dots,7,\,j=0,1\}$, that is all 8 along the line $x_2=1.5$.
The detectors were positioned at $\{(-2+0.4i,1.0+1.0j,0),\
i=0,1,\dots,10,\,j=0,1\}$, that is eleven detectors along each of the two lines
$x_2=1$ and
$x_2=2$.
This corresponds to a mine
search, where the detectors and the sources must be placed outside of
the searched ground. The results of the identification for noise level $\delta=0.00$ in the data
are shown in Figure \ref{subs_fig3} and Table \ref{subs_table4}. The results for noise 
level $\d=0.05$ are shown in Table \ref{subs_table5}.

\begin{table}
\caption{Experiment 2. Identified inclusions, no noise, $\delta=0.00$.}
\begin{tabular}{|r|r|r|r|}
\hline
 $\qquad   x_1$ & $\qquad x_2$ & $\qquad x_3$ & $\qquad v$ \\
\hline
 1.656 & -0.409 & 0.857 & 1.75451 \\
 -1.476 & -0.475 & 0.620 & 0.48823 \\
 1.209 & 0.605 & 0.382 & 0.60886 \\
  -0.225 & 0.469 & 0.266 & 0.69805 \\
-1.406 & 0.228 & 0.159 & 0.59372 \\
\hline
\end{tabular}\label{subs_table4}

\end{table} 

\begin{table}
\caption{Experiment 2. Identified inclusions, $\delta=0.05$.}
\begin{tabular}{|r|r|r|r|}
\hline
 $\qquad   x_1$ & $\qquad x_2$ & $\qquad x_3$ & $\qquad v$ \\
\hline
 1.575 & -0.523 & 0.735 & 1.40827 \\
 -1.628 & -0.447 & 0.229 & 1.46256 \\
 1.197 & 0.785 & 0.578 & 0.53266 \\
 -0.221 & 0.460 & 0.231 & 0.67803 \\
\hline
\end{tabular}\label{subs_table5}

\end{table} 

In general, the execution times were less than 2 minutes on a 333MHz PC.
As it can be seen from the results, the method achieves a perfect
identification in the Experiment \#1 when no noise is present. The
identification deteriorates in the presence of noise, as well as if the
sources and detectors are not located directly above the search area.
Still the inclusions with the highest intensity and the closest ones to
the surface are identified, while the deepest and the weakest are lost.
This can be expected, since their influence on the cost functional is
becoming comparable with the background noise in the data.

In summary, the proposed method for the identification of small inclusions
can be used in geophysics, medicine and technology.
It can be useful in the development of new approaches to ultrasound
mammography. It
can also be used
for localization of holes and cracks in metals and other
materials, as well as for finding mines from
surface measurements of acoustic pressure and possibly in other
problems of interest in various applications.

The HSD minimization method is a specially designed low-dimensional
minimization method, which is well suited for many inverse type
problems. The problems do not necessarily have to be within
the range of applicability of the Born approximation.
It is highly desirable to apply HSD method to practical problems and
to compare its performance with other methods.


\section{Identification of layers in multilayer particles.}

\subsection{Problem Description}
Many practical problems require an identification of the internal
structure of an object given some measurements on its surface. In this
section we study such an identification for a multilayered particle
illuminated by acoustic or electromagnetic plane waves.
Thus the problem discussed here is an inverse scattering problem. A similar
problem for the particle identification from the light scattering data is studied in \cite{zub}.
Our approach is to reduce
the inverse problem to the best fit to data multidimensional
minimization.

Let $D\subset \rc^2$ be the circle of a radius $R>0$,
\begin{equation}\label{21}
D_m=\{x\in \rc^2\, :
r_{m-1}<|x|<r_m\,,\quad m=1,2,\dots ,N\}
\end{equation}
and $S_m=\{x\in \rc^2 : |x|=r_m\}$ for $0=r_0<r_1<\cdots<r_N < R$.
Suppose that a multilayered scatterer in $D$ has a constant refractive
index $n_m$ in the region $D_m\,,\quad m=1,2,\dots ,N $. If the
scatterer is illuminated by a plane harmonic wave
then, after the time dependency is eliminated, the total field
$u(x)=u_0(x)+u_s(x)$  satisfies the Helmholtz equation
\begin{equation}\label{22}
\Delta u+k_0^2u=0\,,\quad |x|>r_N
\end{equation}
where $u_0(x)=e^{ik_0x\cdot\alpha}$ is the incident field and $\alpha$ is
the unit vector in the direction of propagation. The scattered field
$u_s$ is required to satisfy the radiation condition at infinity,
see \cite{r190}.

Let $k_m^2=k_0^2n_m$.
We consider the following transmission problem
\begin{equation}\label{23}
\Delta u_m+k_m^2u_m=0\,\quad x\in D_m\,,
\end{equation}
under the assumption that the fields $u_m$ and their normal derivatives are
continuous across the boundaries $S_m\,, \,m=1,2,\dots ,N$.

In fact, the choice of the boundary conditions on the boundaries $S_m$
depends on the physical model under the consideration. The above model
may or may not be adequate for an electromagnetic or acoustic
scattering, since the model may require additional parameters (such as
the mass density and the compressibility) to be accounted for.
However, the basic computational approach remains the same.
 For more details on transmission problems, including the
questions on the existence and the uniqueness of the solutions,  see
\cite{athrammstrat}, \cite{r410}, and \cite{ewing}.

The Inverse Problem to be solved is:

{\bf IPS:} {\it Given $u(x)$ for all $x\in S=\{x: |x|=R)$
at a fixed $k_0>0$, find
the number $N$ of the layers, the location of the layers, and their refractive indices
$n_m \,, \,m=1,2,\dots ,N$ in \eqref{23}.
}

Here the IPS stands  for a Single frequency Inverse Problem.
 Numerical experience shows that there
are some practical difficulties in the successful resolution of the IPS
even when no noise is present, see \cite{gut5}.
While there are some results on the uniqueness for the IPS (see
\cite{athrammstrat,r410}), assuming that the refractive indices are known, and only
the layers are to be identified, the stability estimates are few,
\cite{r325},\cite{r329}, \cite{r425}.
The identification is successful, however, if the scatterer is subjected
to a probe with plane waves of several frequencies. Thus we state the Multifrequency
Inverse Problem:

{\bf IPM:} {\it Given $u^p(x)$ for all $x\in S=\{x: |x|=R)$
at a finite number $P$ of wave numbers $k_0^{(p)}>0$, find
the number $N$ of the layers, the location of the layers, and their refractive indices
$n_m \,, \,m=1,2,\dots ,N$ in \eqref{23}.
}

\subsection{Best Fit Profiles and Local Minimization Methods}
If the refractive indices  $n_m$ are sufficiently close to $1$, then we say that the
scattering is weak. In this case the scattering is
described by the Born
approximation, and there are methods for the solution
of the above Inverse Problems. See \cite{coltonmonk},
\cite{r190} and \cite{r313} for further
details. In particular,  the Born inversion is
an ill-posed problem even if the Born approximation is very accurate,
see \cite{rammb2},
 or \cite{r249}. When the assumption of the Born approximation is not
appropriate, one
matches the given observations to a set of
solutions for the Direct Problem. Since our interest is in the solution of
the IPS and IPM in the non-Born region of scattering, we choose to
follow the best fit to data approach. This approach is used widely in a
variety of applied problems, see e. g. \cite{biegler}.

Note, that, by the assumption, the scatterer has the rotational
symmetry. Thus we only need to know the data for one direction of the incident
plane wave. For this reason we fix $\alpha=(1,0)$ in \eqref{22} and
define
the (complex) functions
\begin{equation}
g^{(p)}(\theta)\,,\quad 0\leq \theta< 2\pi, \quad p=1,2,\dots ,P,
\end{equation}
to be the observations
measured on the surface $S$ of the ball $D$ for a finite set of free space
wave numbers $k_0^{(p)}$.

Fix a positive integer $M$. Given a configuration
\begin{equation}
Q=(r_1,r_2,\dots,r_M,n_1,n_2,\dots,n_M)
\end{equation}
we solve the Direct Problem \eqref{22}-\eqref{23} (for each free space wave number  $k_0^{(p)}$)
 with the layers
$D_m=\{x\in \rc^2\, :
r_{m-1}<|x|<r_m\,,\quad m=1,2,\dots ,M\}$, and the corresponding refractive indices
$n_m$, where $r_0=0$.
Let
\begin{equation}
w^{(p)}(\theta)=u^{(p)}(x)\big|_{x\in S}\,.
\end{equation}

Fix a set of angles $\Theta=(\theta_1,\theta_2,\dots,\theta_L)$
and  let
\begin{equation}
\|w\|_2=(\sum_{l=1}^L w^2(\theta_l))^{1/2}.
\end{equation}

Define
\begin{equation}\label{28}
\Phi(r_1,r_2,\dots,r_M,n_1,n_2,\dots,n_M)=\frac{1}{P}\sum^P_{p=1}\frac{\|w^{(p)}-
g^{(p)}\|_2^2}{\|g^{(p)}\|_2^2}\,,
\end{equation}
where the same set $\Theta$ is used for $g^{(p)}$ as for $w^{(p)}$.

We solve the IPM by minimizing the above best fit to data functional $\Phi$
over an appropriate set of admissible parameters $A_{adm}\subset
\rc^{2M}$.

It is reasonable to assume that the underlying physical problem gives
some estimate for the bounds $n_{low}$ and $n_{high}$
of the refractive indices $n_m$ as well as for the bound $M$ of
the expected number of layers $N$. Thus,

\begin{equation}
A_{adm} \subset \{(r_1,r_2,\dots,r_M,n_1,n_2,\dots,n_M)\ : \ 0\leq r_i\leq R\,,\
n_{low}\leq n_m \leq n_{high}\}\,.
\end{equation}

Note, that the admissible configurations must also satisfy

\begin{equation}
r_1\leq r_2\leq r_3 \leq\dots\leq r_M\,.
\end{equation}

It is well known that  a
multidimensional minimization is a difficult problem,
unless the objective function is "well behaved". The most important
quality of such a cooperative function is the presence of just a few
local minima.
Unfortunately, this is, decidedly, not the case in many
applied problems, and, in particular, for the problem under the consideration.

To illustrate this point further, let  $P$ be the set of three free space wave
numbers $k_0^{(p)}$  chosen to be

\begin{equation}
P=\{3.0,\ 6.5,\ 10.0\}\,.
\end{equation}

\begin{figure}[tb]
\vspace{5pc}
\includegraphics*{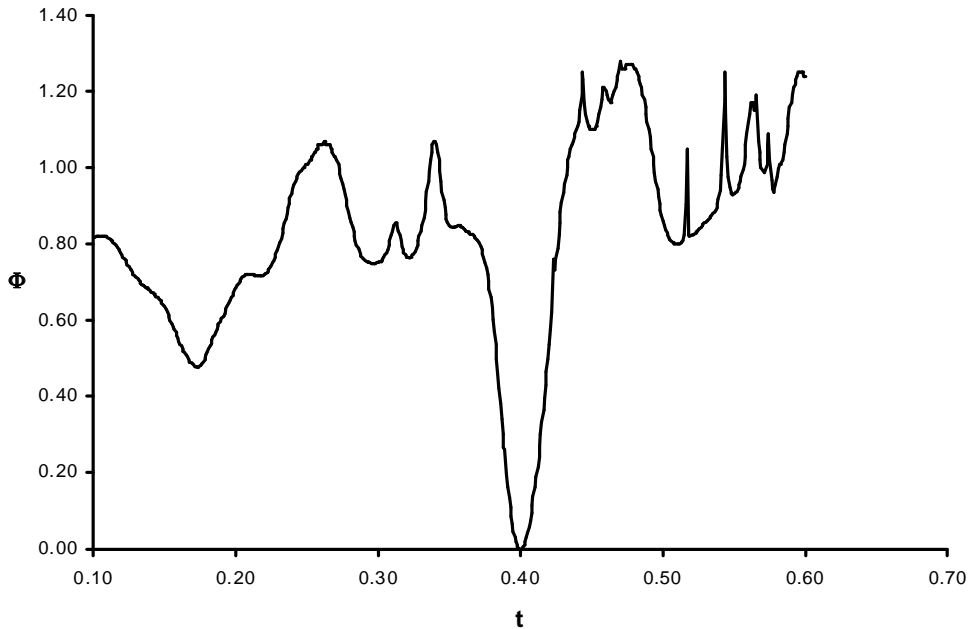}
\caption{Best fit profile for the configurations $q_t$; Multiple frequencies
$P=\{3.0,\ 6.5,\ 10.0\}$.} \label{g5fig5}
\end{figure}

Figure \ref{g5fig5} shows the profile of the functional $\Phi$ as a function of the
variable $t\,,\, 0.1\leq t \leq 0.6$ in the configurations $q_t$ with

\[
n(x)=\begin{cases}
0.49 & 0\leq |x| < t\\
9.0 & t \leq |x| < 0.6\\
1.0 & 0.6\leq |x| \leq 1.0
\end{cases}
\]

Thus the objective function $\Phi$ has many
local minima even along this arbitrarily chosen one dimensional
cross-section of the admissible set. There are sharp peaks and large
gradients.
Consequently, the gradient based methods (see \cite{bre,des,fletcher,hestenes,jac,pol}),
 would not be successful
for a significant portion of this region. It is also appropriate to
notice that the dependency of $\Phi$ on its arguments is highly
nonlinear. Thus, the gradient computations have to be done numerically, which
makes them computationally expensive. More
importantly, the gradient based minimization methods (as expected)
perform poorly for these problems.

These complications are avoided by considering conjugate gradient type algorithms
 which do not require the
knowledge of the derivatives at all, for example the Powell's method.
Further refinements in the deterministic phase of the minimization algorithm
are needed to achieve more consistent performance. They include special
line minimization, and Reduction procedures similar to the ones
discussed in a previous section on the identification of underground
inclusions. We skip the details and refer the reader to \cite{gut5}.

In summary, the entire Local Minimization Method {\bf (LMM)} consists of the
following:

\subsubsection*{Local Minimization Method (LMM)}
\begin{enumerate}

\item  Let your starting configuration be
$Q_0=(r_1,r_2,\dots,r_M,n_1,n_2,\dots,n_M)$.

\item  Apply the Reduction Procedure to $Q_0$, and obtain a reduced configuration
$Q_0^r$ containing $M^r$ layers.

\item  Apply the Basic Minimization Method in $A_{adm}\bigcap \rc^{2M^r}$
with the starting point $Q_0^r$, and obtain a configuration $Q_1$.

\item  Apply the Reduction Procedure to $Q_1$, and obtain a final reduced configuration
$Q_1^r$.

\end{enumerate}

\subsection{Global Minimization Methods}

Given an initial configuration $Q_0$ a local minimization method finds a
local minimum near $Q_0$. On the other hand, global minimization methods
explore the entire admissible set to find a global minimum of the
objective function. While the local minimization is, usually,
deterministic, the majority of the global methods are probabilistic in
their nature. There is a great interest and activity in the development of efficient
global minimization methods, see e.g. \cite{biegler},\cite{bomze}. Among them are the
simulated annealing method (\cite{kgv},\cite{kir}), various genetic algorithms \cite{hh},
 interval method,
TRUST method (\cite{bar},\cite{bpr}), etc. As we have already mentioned before, the best fit to
data functional $\Phi$ has many narrow local minima. In this situation
it is exceedingly unlikely to get the minima points by chance alone.
Thus our special interest is for the minimization methods, which combine
a global search with a local minimization. In \cite{r399} we
developed such a method (the Hybrid Stochastic-Deterministic Method),
and applied it for the identification of small subsurface particles,
provided a set of surface measurements, see Sections 2.2-2.4. The HSD method could be
classified as a variation of a genetic algorithm with a local search
with reduction.
In this paper we consider the performance of two
algorithms: Deep's Method, and Rinnooy Kan and Timmer's
Multilevel Single-Linkage Method. Both combine a global and a local search
to determine a global minimum. Recently these methods have been applied
to a similar problem of
the identification of  particles from their light scattering characteristics in \cite{zub}.
Unlike \cite{zub}, our experience shows that
Deep's method has failed consistently for the type of problems we are
considering. See \cite{Deep} and \cite{zub} for more details on Deep's
Method.

\subsubsection*{Multilevel Single-Linkage Method (MSLM)}

Rinnooy Kan and Timmer in \cite{rkt1} and \cite{rkt2} give a detailed
description of this algorithm. Zakovic et. al. in \cite{zub} describe in detail an
experience of its application to an inverse light scattering problem.
They also discuss different stopping criteria for the MSLM.
Thus, we only give here a shortened and an informal description of this method and of its
algorithm.

In a pure {\bf Random Search} method a batch $H$ of $L$ trial points is generated in
$A_{adm}$ using a uniformly distributed random variable. Then a local
search is started from each of these $L$ points. A local minimum with
the smallest value of $\Phi$ is declared to be the global one.

A refinement of the Random Search is the {\bf Reduced Sample Random Search} method.
Here we use
only a certain fixed fraction $\gamma<1$ of the original batch of $L$
points to proceed with the local searches. This reduced sample $H_{red}$ of
$\gamma L$ points is chosen to contain the points with the smallest
$\gamma L$ values of $\Phi$ among the original batch. The local searches
are started from the points in this reduced sample.

Since the local searches dominate the computational costs, we would like
to initiate them only when it is truly necessary. Given a critical
distance $d$ we define a cluster to be a group of points located within
the distance $d$ of each other. Intuitively, a local search started
from the points within a cluster should result in the same local
minimum, and, therefore, should be initiated only once in each cluster.

Having tried all the points in the reduced sample we have an information
on the number of local searches performed and the number of local minima
found. This information and the critical distance $d$ can be used to
determine a statistical level of confidence, that all the local minima
have been found. The algorithm is terminated (a stopping criterion is satisfied)
if an a priori level of confidence is reached.

If, however, the stopping criterion is not satisfied, we perform another
iteration of the MSLM by generating another batch of $L$ trial points.
Then it is
combined with the previously generated batches to obtain an enlarged batch $H^j$ of
$jL$ points (at iteration $j$), which leads to a
reduced sample $H^j_{red}$ of $\gamma jL$ points. According to MSLM the
critical distance $d$ is reduced to $d_j$, (note that $d_j\ra 0$ as $j\ra\infty$, since we want to
find a minimizer), a local minimization is attempted once within each cluster,
the information on the number of local minimizations performed and the
local minima found is used to determine if the algorithm should be
terminated, etc.

The following is an adaptation of the MSLM method to the inverse scattering problem
presented in Section 3.1. The LMM local
minimization method introduced in the previous Section is used here to
perform local searches.

\subsubsection*{MSLM}
(at iteration $j$).
\begin{enumerate}

\item  Generate another batch of $L$ trial points (configurations) from a random
uniform distribution in $A_{adm}$.
Combine it with the previously generated batches to obtain an enlarged batch $H^j$ of
$jL$ points.

\item Reduce $H^j$ to the reduced sample $H^j_{red}$ of $\gamma jL$
points, by selecting the points with the smallest $\gamma jL$ values of $\Phi$
in $H^j$.

\item Calculate the critical distance $d_j$ by

\[
d_j^r=\pi^{-1/2}\left( \Gamma\left(1+\frac{M}{2}\right)R^M\frac{\sigma\ln jL}
{jL}\right)^{1/M}\,,
\]
\[
d_j^m=\pi^{-1/2}\left( \Gamma\left(1+\frac{M}{2}\right)(n_{high}-n_{low})^M\frac{\sigma\ln jL}
{jL}\right)^{1/M}\,.
\]
\[
d_j=\sqrt{(d_j^r)^2+(d_j^n)^2}
\]

\item Order the sample points in $H^j_{red}$ so that
$\Phi(Q_i)\leq\Phi(Q_{i+1})$, $i=1,\dots,\gamma jL$.
For each value of $i$, start the local minimization from $Q_i$, unless
there exists an index $k<i$, such that $\|Q_k-Q_i\|\leq d_j$. Ascertain if the
result is a known local minimum.

\item Let $K$ be the number of local minimizations performed, and $W$ be
the number of different local minima found. Let
\[
W_{tot}=\frac{W(K-1)}{K-W-2}
\]
The algorithm is terminated if

\begin{equation}\label{29}
W_{tot}<W+0.5\,.
\end{equation}

\end{enumerate}

Here $\Gamma$ is the gamma function, and $\sigma$ is a fixed constant.

A related algorithm (the Mode Analysis) is based on a subdivision of the
admissible set into smaller volumes associated with local minima. This
algorithm is also discussed in \cite{rkt1} and \cite{rkt2}. From the
numerical studies presented there, the authors deduce their preference
for the MSLM.

The presented MSLM algorithm was successful in the identification of various 2D
layered particles, see \cite{gut5} for details.


\section{Potential scattering and the Stability Index method.}

\subsection{Problem description} Potential scattering
problems are important in quantum mechanics, where they
appear in the context of scattering of particles bombarding
an atom nucleus. One is interested in reconstructing the
scattering potential from the results of a scattering
experiment. The examples in  Section 4 deal with finding a 
spherically symmetric
($q=q(r),\; r=|x|$) potential from the 
fixed-energy scattering data, which in this case consist 
of the fixed-energy phase shifts. In \cite{R8}, \cite{r425}, 
\cite{r470} and \cite{r474}
 the three-dimensional inverse scattering problem with 
fixed-energy data is treated.

Let $q(x),\, x\in \rc^3,$ be a real-valued potential with compact
support. Let $R>0$ be a number
 such that $q(x)=0$ for $\abs{x} > R$. We also assume that
$q\in L^2(B_R)\,,\ B_R=\{x:\abs{x}\leq R, x\in \rc^3\}$. Let
$S^2$ be the unit sphere, and $\alpha \in S^2$. For a given
energy $k>0$ the scattering solution $\psi(x,\alpha)$ is
defined as the solution of

\begin{equation}
\Delta \psi+k^2\psi-q(x)\psi= 0\,,\quad x \in \rc^3
\end{equation}
satisfying the following asymptotic condition at infinity:

\begin{equation}
\psi=\psi_0+v,\quad \psi_0:=e^{ik\alpha\cdot x}\,,\quad \alpha\in S^2\,,
\end{equation}

\begin{equation}
\lim_{r\rightarrow\infty}\int_{\abs{x}=r}\left| \frac{\partial v}{\partial
r}-ikv\right|^2 ds=0\,.
\end{equation}

It can be shown, that

\begin{equation}
\psi(x,\alpha)=\psi_0+A(\alpha',\alpha,k)
\frac{e^{ikr}}{r}+o\left(\frac{1}{r}\right)\,,\;
\text{as}\ \ r\rightarrow\infty\,, \quad
\frac{x}{r}=\alpha'\, \quad r:=|x|.
\end{equation}

The function $A(\alpha',\alpha,k)$ is called
the scattering amplitude, $\alpha$ and $\alpha'$
are the directions of the incident and scattered waves, and $k^2$
is the energy, see  \cite{n}, \cite{r313}.

For spherically symmetric scatterers $q(x)=q(r)$ the scattering
amplitude satisfies $A(\alpha',\alpha,k)=A(\alpha'\cdot\alpha,k)$.
The converse is established in
\cite{r6}. Following \cite{r394}, the scattering amplitude
for $q=q(r)$ can be written as

\begin{equation}
A(\alpha',\alpha,k)=\sum^\infty_{l=0}\sum^l_{m=-
l}A_l(k)Y_{lm}(\alpha')\overline{Y_{lm}(\alpha)}\,,
\end{equation}
where $Y_{lm}$ are the spherical harmonics, normalized in
$L^2(S^2),$ and the bar denotes the
complex conjugate.

The fixed-energy phase shifts $-\pi<\delta_l\leq\pi$
($\delta_l=\delta(l,k)$,  $k>0$ is fixed) are related to $A_l(k)$
(see e.g., \cite{r394}) by the formula:

\begin{equation}
A_l(k)=\frac{4\pi}{k}e^{i\delta_l}\sin(\delta_l)\,.
\end{equation}

Several parameter-fitting procedures were proposed for calculating the
potentials from the fixed-energy phase shifts, (by Fiedledey,
Lipperheide, Hooshyar and Razavy, Ioannides and Mackintosh, Newton,
 Sabatier, May and Scheid, Ramm and others).
These works  are referenced and their results are described in \cite{cs}
and \cite{n}. Recent works \cite{gut3,gut5,r399,r423} and \cite{gutmanramm2,
r394,rsmi}, \cite{r394} present
new numerical methods for solving this problem. In \cite{r431} (also see  \cite{r460,r470})
it is proved that the R.Newton-P.Sabatier method for solving inverse scattering problem
the fixed-energy phase shifts as the data (see \cite{n,cs} )
is fundamentally wrong in the sense that its foundation is wrong.
In  \cite{R5} a counterexample is given to a
uniqueness theorem claimed in a modification of the R.Newton's
inversion scheme.

Phase shifts for a spherically symmetric potential can be
computed by a variety of methods, e.g. by a variable phase method
described in \cite{calogero}. The computation involves solving
a nonlinear ODE for each phase shift. However, if the potential is compactly
supported and piecewise-constant, then
a much simpler method described in \cite{ars} and \cite{r435} can be
used. We refer the reader to these papers for details.

Let $q_0(r)$ be a spherically symmetric piecewise-constant potential,
 $\{\tilde\delta(k,l)\}_{l=1}^N$ be the set of its phase shifts
for a fixed $k>0$ and a sufficiently large $N$.
Let $q(r)$ be another  potential, and let
 $\{\delta(k,l)\}_{l=1}^N$ be the set of its phase shifts.

The best fit to data
function $\Phi(q,k)$ is defined by

\begin{equation}\label{phi}
\Phi(q,k)=\frac{\sum^N_{l=1}\abs{\delta(k,l)-
\tilde\delta(k,l)}^2}
{\sum^N_{l=1}\abs{\tilde\delta(k,l)}^2}.
\end{equation}

The phase shifts are known to decay rapidly with $l$, see \cite{rai}.
Thus, for sufficiently large $N$, the function $\Phi$ is practically the same as
the one which would use all the shifts in (\ref{phi}).
The inverse problem of the reconstruction of the potential from its
fixed-energy phase shifts is reduced to the minimization of the
objective function $\Phi$
over an appropriate admissible set.

\subsection{Stability Index Minimization Method}

Let the minimization problem be
\begin{equation}\label{sta1}
\min\{\Phi(q) \ : \ q\in A_{adm}\}
\end{equation}

Let $\tilde q_0$ be its global minimizer. Typically, the
structure of the objective function $\Phi$ is quite complicated: this
function may have
many local minima. Moreover, the objective function in a neighborhood of
 minima can be nearly
flat resulting in large minimizing sets defined by

\begin{equation}\label{sta2}
S_{\ep}=\{q\in A_{adm}\ :\ \Phi(q)<\Phi(\tilde q_0)+\ep \}
\end{equation}
for an $\ep>0$.

{\it Given an $\ep>0$, let $D_\ep$ be the diameter of the minimizing set
$S_\ep$, which we call the {\bf Stability Index $D_\ep$} of the 
minimization
problem (\ref{sta1}).}

Its usage is explained below.

One would expect to obtain stable identification for minimization
problems with small (relative to the admissible set) stability indices. Minimization
problems with large stability indices have distinct minimizers with
practically the same values of the objective function. If no additional
information is known, one has an uncertainty
of the minimizer's choice. The stability index provides a quantitative
measure of this uncertainty or instability of the minimization.

If  $D_\ep<\eta$, where $\eta$ is an a priori chosen treshold,
then one can solve  the global minimization problem stably.
In the above general scheme it is not discussed in detail what
are possible algorithms for computing the Stability Index.

One idea to construct such an algorithm is to
iteratively estimate stability indices of the minimization
problem, and, based on this information, to conclude if the method has
achieved a stable minimum.

One such algorithm is an Iterative Reduced Random Search (IRRS) method,
which uses the Stability Index for its stopping criterion.
Let a batch $H$ of $L$ trial points be randomly generated in the admissible set
$A_{adm}$.
Let   $\gamma$ be a certain fixed fraction, e.g., $\gamma=0.01$.
 Let
$S_{min}$ be the subset of $H$ containing points
 $\{p_i\}$ with the smallest $\gamma L$ values of the objective function
$\Phi$ in $H$. We call $S_{min}$ the minimizing set.  If all the
minimizers in $S_{min}$ are close to each other, then the objective
function $\Phi$ is not flat near the global minimum. That is, the method
identifies the minimum consistently.
 Let
$\|\cdot\|$ be a norm in the admissible set.

Let
\[\ep=\max_{p_j \in S_{min}}\Phi(p_j)-\min_{p_j \in S_{min}}\Phi(p_j)
\]
and
\begin{equation}\label{diam}
\tilde D_\ep=diam(S_{min})=\max\{\|p_i-p_j\|\ :\ p_i,p_j\in S_{min}\}\,.
\end{equation}

Then $\tilde D_\ep$ can be considered an estimate for the {\bf Stability Index} 
$D_\ep$ of
the minimization problem. The Stability Index reflects the size of the minimizing sets.
 Accordingly, it is used as a self-contained stopping criterion for an
iterative minimization procedure. 
The identification is considered to be stable
if the Stability Index $D_\ep<\eta$, for an a priori chosen
$\eta>0$. Otherwise, another batch of $L$ trial
points is generated, and the process is repeated. We used $\beta=1.1$ as
described below in the stopping criterion to determine if subsequent
iterations do not produce a meaningful reduction of the objective
function.

More precisely

\subsection*{Iterative Reduced Random Search (IRRS)} (at the $j-$th
iteration).

Fix $0<\gamma <1,\ \beta > 1,\ \eta>0$ and $N_{max}$.

\begin{enumerate}
\item  Generate another batch $H^j$ of $L$ trial points in $A_{adm}$ using a
random distribution.

\item Reduce $H^j$ to the reduced sample $H^j_{min}$ of $\gamma L$ points
by selecting the points in $H^j$ with the smallest $\gamma L$ values of
$\Phi$.

\item Combine $H^j_{min}$ with $H^{j-1}_{min}$ obtained at the previous
iteration. Let $S^j_{min}$ be the set of $\gamma L$
points from $H^j_{min}\cup H^{j-1}_{min}$ with the smallest
 values of $\Phi$. (Use $H^1_{min}$ for $j=1$).

\item Compute the Stability Index (diameter) $D^j$ of $S^j_{min}$ by
$D^j=\max\{\|p_i-p_k\|\ :\ p_i,p_k\in S_{min}\}\,.$

\item Stopping criterion.

Let $p\in S^j_{min}$ be the point with the smallest value
of $\Phi$ in $S^j_{min}$ (the global minimizer).

If $D^j\leq\eta$, then stop. The global minimizer is $p$. The
minimization is stable.

If $D^j>\eta$ and $\Phi(q)\leq\beta \Phi(p)\ : \ q\in S^j_{min}$, then
stop. The minimization is unstable. The Stability Index $D^j$ is the
measure of the instability of the minimization.

Otherwise, return to step 1
and do another iteration, unless the maximum number of iterations $N_{max}$ is
exceeded.

\end{enumerate}

One can make the stopping criterion more meaningful by computing a
normalized stability index. This can be achieved by dividing $D^j$ by a
fixed normalization constant, such as the diameter of the entire
admissible set $A_{adm}$.
To improve the performance of the algorithm in specific problems we
found it useful to modify (IRRS) by combining the stochastic (global)
search with a deterministic local minimization. Such
 Hybrid Stochastic-Deterministic (HSD) approach
 has proved to be successful for a variety of problems
in inverse quantum scattering (see \cite{r435, gut5,gutmanramm2})
as well as in other
applications (see \cite{r399,gut3}). A somewhat different
implementation of the Stability Index Method is described in
\cite{r423}.

We seek the potentials $q(r)$ in the class of piecewise-constant, spherically
symmetric real-valued functions. Let the admissible set be
\begin{equation}\label{adm}
A_{adm} \subset \{(r_1,r_2,\dots,r_M,q_1,q_2,\dots,q_M)\ : \ 0\leq r_i\leq R\,,\
q_{low}\leq q_m \leq q_{high}\}\,,
\end{equation}

where the bounds $q_{low}$ and $q_{high}$
for the potentials, as well as the bound $M$ on
the expected number of layers are assumed to be known.

A configuration $(r_1,r_2,\dots,r_M,q_1,q_2,\dots,q_M)$ corresponds to
the potential

\begin{equation}
q(r)=q_m\,,\quad \text{for}\quad r_{m-1}\leq r<r_m\,,\quad 1\leq m\leq M\,,
\end{equation}
where $r_0=0$ and $q(r)=0$ for $r\geq r_M=R$.

Note, that the admissible configurations must also satisfy
\begin{equation}\label{admr}
r_1\leq r_2\leq r_3 \leq\dots\leq r_M\,.
\end{equation}

We used $\beta=1.1$, $\epsilon=0.02$ and $j_{max}=30$. The choice of
these and other parameters
 ($L=5000,\, \gamma=0.01,\ \nu=0.16$ )
 is dictated by their meaning in the algorithm and the comparative performance
of the program at their different values. As usual,
some adjustment of the parameters, stopping criteria,
etc., is needed to achieve the optimal performance of the algorithm.
The deterministic part of the IRRs algorithm was based on the Powell's
minimization method, one-dimensional minimization, and a Reduction
procedure similar to ones described in the previous section 3, see
\cite{r435} for details.

\subsection{Numerical Results}
We studied the performance of the algorithm for 3
different potentials $q_i(r),\ i=1,2,3$ chosen from the physical
considerations.

The potential $q_3(r) = - 10$ for $0\leq r < 8.0$ and $q_3 = 0 $
for $r\geq 8.0$ and a wave number $k=1$ constitute a typical
example  for elastic scattering of neutral particles
in nuclear and atomic physics. In nuclear physics one measures the length in
units of fm = $10^{-15}$m, the quantity $q_3$ in units of
1/fm$^2$, and the wave number in units of 1/fm. The physical
potential and incident energy are given by $V(r) = \frac
{\hbar^2}{2\mu} q_3(r)$ and $E = \frac {\hbar^2 k^2}{2\mu}$,
respectively.
here $\hbar:= \frac {h}{2\pi}$,  $h=6.625 10^{-27}$ erg$\cdot$s
is the Planck constant,
 $\hbar c=197.32$  MeV$\cdot$fm, $c=3\cdot 10^6$ m/sec is the
velocity of light,
and $\mu$ is the mass of a neutron.
By choosing the mass $\mu$ to be equal to
the mass of a neutron $\mu$ = 939.6 MeV/$c^2$, the potential and
energy have the values of $V(r)$ = -207.2 MeV for $0\leq r < 8.0$ fm
and $E( k=$1/fm ) = 20.72 MeV. In atomic physics one uses atomic
units with the Bohr radius $a_0 = 0.529\cdot10^{-10}$m as the unit of
length. Here, $r, k $ and $q_3$ are measured in units of $a_0,
1/a_0$ and $1/a_0^2$, respectively. By assuming a scattering of an
electron with mass $m_0$ = 0.511 MeV/$c^2$, we obtain the potential
and energy as follows: $V(r) = -136$ eV for $0\leq r < 8 a_0 =
4.23\cdot10^{-10}$m and $E( k=1/a_0 ) = 13.6$ eV. These numbers
give motivation for the choice of examples applicable in
nuclear and atomic physics.

The method used here deals with finite-range
(compactly supported) potentials. One can use this method for
potentials with the Coulomb tail or other potentials
of interest in physics, which are not of finite range.
This is done by using the phase shifts transformation
method which allows one to transform the phase shifts corresponding to
a potential, not of finite range, whose behavior is known for
$r>a$, where $a$ is some radius, into the phase shifts corresponding
to a potential of finite range $a$ (see \cite{apagyi}, p.156).

In practice differential cross section is measured
at various angles, and from it the fixed-energy phase
 shifts are calculated by a parameter-fitting procedure.
Therefore, we plan in the future work to generalize
the stability index method to the
case when the original data are the values
of the differential cross section, rather than
the phase shifts.



For the physical reasons discussed above,
we choose the following three potentials:
\[
q_1(r)=\begin{cases}
-2/3 & 0\leq r < 8.0\\
0.0 &  r \geq 8.0
\end{cases}
\]

\[
q_2(r)=\begin{cases}
-4.0 & 0\leq r < 8.0\\
0.0 &  r \geq 8.0
\end{cases}
\]

\[
q_3(r)=\begin{cases}
-10.0 & 0\leq r < 8.0\\
0.0 &  r \geq 8.0
\end{cases}
\]

In each case the following values of the parameters have been used. The
radius $R$ of the support of each $q_i$ was chosen to be $R=10.0$.
The admissible set $A_{adm}$ (\ref{adm}) was defined with $M=2$. The
Reduced Random Search parameters:
$L=5000\,,\;\gamma=0.01\,,\;\nu=0.16\,,\; \epsilon=0.02\,,\;\beta=1.10\,,j_{max}=30$.
The value $\epsilon_r=0.1$ was used in
the Reduction Procedure during the local minimization phase.
The initial configurations were generated using a random
number generator with seeds determined by the system time.
A typical run time was about 10 minutes on a 333 MHz PC,
depending on the number of iterations in IRRS.
The number $N$ of the shifts used in (\ref{phi}) for the formation of
the objective function $\Phi(q)$ was $31$ for all the wave numbers.
It can be seen that the shifts for the potential $q_3$ decay rapidly for $k=1$,
but they remain large for $k=4$.
 The upper and lower bounds for the potentials $q_{low}=-20.0$ and $q_{high}=0.0$
used in the definition of the admissible set $A_{adm}$ were chosen to
reflect a priori information about the potentials.


The identification was attempted with 3 different noise levels $h$.
The levels are $h=0.00$ (no noise), $h=0.01$ and $h=0.1$.
 More precisely, the noisy phase shifts $\delta_h(k,l)$ were obtained from
the exact phase shifts $\delta(k,l)$ by the formula

\[
\delta_h(k,l)=\delta(k,l)(1+(0.5-z)\cdot h)\,,
\]
where
$z$ is the uniformly distributed on $[0,1]$ random variable.

The distance $d(p_1(r),p_2(r))$ for potentials in step 5 of the IRRS algorithm was
computed as

\[
 d(p_1(r),p_2(r))=\|p_1(r)-p_2(r)\|\,
\]
where the norm is the $L_2$-norm in $\rc^3$.

The results of the identification algorithm (the Stability Indices) for
different iterations of the IRRS algorithm
are shown in Tables \ref{table_stab2}-\ref{table_stab4}.


\begin{table}
\caption{Stability Indices for $q_1(r)$ identification at different noise levels $h$.}

\begin{tabular}{r r r r r}
\hline

$k$ & $Iteration$ & $h=0.00$ & $h=0.01$ & $h=0.10$\\

\hline

\hline

1.00 & 1 & 1.256985 & 0.592597   & 1.953778  \\
   & 2 &    0.538440 & 0.133685 & 0.799142 \\
   & 3 &    0.538253 & 0.007360 & 0.596742 \\
   & 4 &    0.014616 & & 0.123247 \\
 &5 &&& 0.015899 \\
\hline
2.00 & 1 & 0.000000 & 0.020204   & 0.009607  \\

\hline
2.50 & 1 & 0.000000 & 0.014553   & 0.046275  \\
\hline
3.00 & 1 & 0.000000 & 0.000501   & 0.096444  \\
\hline
4.00 & 1 & 0.000000 & 0.022935  & 0.027214   \\
\hline

\hline
\end{tabular}\label{table_stab2}
\end{table}


\begin{table}
\caption{Stability Indices for $q_2(r)$ identification at different noise levels $h$.}

\begin{tabular}{r r r r r}
\hline

$k$ & $Iteration$ & $h=0.00$ & $h=0.01$ & $h=0.10$\\

\hline

\hline

1.00 & 1 & 0.774376 & 0.598471   & 0.108902 \\
   & 2 &    0.773718 & 1.027345 & 0.023206 \\
   & 3 &    0.026492 & 0.025593 & 0.023206 \\
   & 4 &    0.020522 & 0.029533 & 0.024081 \\
   & 5 &    0.020524 & 0.029533 & 0.024081 \\
   & 6 &    0.000745 & 0.029533 &  \\
   & 7 &  & 0.029533 &  \\
   & 8 &  & 0.029533 &  \\
   & 9 &  & 0.029533 &  \\
   & 10 &  & 0.029533 &  \\
   & 11 &  & 0.029619 &  \\
   & 12 &  & 0.025816 &  \\
   & 13 &  & 0.025816 &  \\
   & 14 &  & 0.008901 &  \\

\hline
2.00 & 1 & 0.863796 & 0.799356   & 0.981239 \\
 & 2 & 0.861842 & 0.799356 & 0.029445 \\
 & 3 & 0.008653 & 0.000993 & 0.029445 \\
 & 4 &  &  & 0.029445 \\
 & 5 &  &  & 0.026513 \\
 & 6 &  &  & 0.026513 \\
 & 7 &  &  & 0.024881 \\

\hline

2.50 & 1 &  1.848910 & 1.632298  & 0.894087 \\
 & 2 & 1.197131 & 1.632298 & 0.507953 \\
 & 3 & 0.580361 & 1.183455 & 0.025454 \\
 & 4 & 0.030516 & 0.528979 & \\
 & 5 & 0.016195 & 0.032661 & \\

\hline
3.00 & 1 & 1.844702 & 1.849016   & 1.708201 \\
 & 2 & 1.649700 & 1.782775   & 1.512821 \\
 & 3 & 1.456026 & 1.782775   & 1.412345 \\
 & 4 & 1.410253 & 1.457020   & 1.156964 \\
 & 5 & 0.624358 & 0.961263   & 1.156964 \\
 & 6 & 0.692080 & 0.961263   & 0.902681 \\
 & 7 & 0.692080 & 0.961263   & 0.902681 \\
 & 8 & 0.345804 & 0.291611   & 0.902474 \\
 & 9 & 0.345804 & 0.286390   & 0.159221 \\
 & 10 & 0.345804 & 0.260693  & 0.154829 \\
 & 11 & 0.043845 & 0.260693  & 0.154829 \\
 & 12 & 0.043845 & 0.260693  & 0.135537 \\
 & 13 & 0.043845 & 0.260693  & 0.135537 \\
 & 14 & 0.043845 & 0.260693  & 0.135537 \\
 & 15 & 0.042080 & 0.157024  & 0.107548 \\
 & 16 & 0.042080 & 0.157024  &  \\
 & 17 & 0.042080 & 0.157024  &  \\
 & 18 & 0.000429 & 0.157024  &  \\
 & 19 & & 0.022988 &  \\

\hline
4.00 & 1 & 0.000000 & 0.000674  & 0.050705 \\
\hline

\end{tabular}\label{table_stab3}
\end{table}


\begin{table}
\caption{Stability Indices for $q_3(r)$ identification at different noise levels $h$.}

\begin{tabular}{r r r r r}
\hline

$k$ & $Iteration$ & $h=0.00$ & $h=0.01$ & $h=0.10$\\

\hline

\hline

1.00 & 1 & 0.564168 & 0.594314   & 0.764340 \\
   & 2 &   0.024441  & 0.028558 & 0.081888 \\
   & 3 &   0.024441  & 0.014468 & 0.050755 \\
   & 4 &   0.024684  &  &  \\
   & 5 &   0.024684  &  &  \\
   & 6 &  0.005800   &  &  \\

\hline

2.00 & 1 & 0.684053 & 1.450148   & 0.485783 \\
   & 2 &  0.423283   & 0.792431 & 0.078716 \\
   & 3 &  0.006291   & 0.457650 & 0.078716 \\
   & 4 &     & 0.023157 & 0.078716 \\
   & 5 &     &  & 0.078716 \\
   & 6 &     &  & 0.078716 \\
   & 7 &     &  & 0.078716 \\
   & 8 &     &  & 0.078716 \\
   & 9 &     &  & 0.078716 \\
   & 10 &     &  & 0.078716 \\
   & 11 &     &  & 0.078716 \\

\hline

2.50 & 1 & 0.126528 & 0.993192   & 0.996519 \\
   & 2 &  0.013621   & 0.105537 & 0.855049 \\
   & 3 &     & 0.033694 & 0.849123 \\
   & 4 &     & 0.026811 & 0.079241 \\
\hline

3.00 & 1 & 0.962483 & 1.541714   & 0.731315 \\
   & 2 &  0.222880   & 0.164744 & 0.731315 \\
   & 3 &  0.158809   & 0.021775 & 0.072009 \\
   & 4 &  0.021366   &  &  \\
   & 5 &  0.021366   &  &  \\
   & 6 &  0.001416   &  &  \\

\hline
4.00 & 1 & 1.714951 & 1.413549   & 0.788434 \\
   & 2 &  0.033024   & 0.075503 & 0.024482 \\
   & 3 &  0.018250   & 0.029385 &  \\
   & 4 &     & 0.029421 &  \\
   & 5 &     & 0.029421 &  \\
   & 6 &     & 0.015946 &  \\

\hline

\end{tabular}\label{table_stab4}
\end{table}

For example, Table \ref{table_stab4} shows that for $k=2.5,\ h=0.00$ the Stability Index has reached
the value $0.013621$ after 2 iteration. According to the Stopping
criterion for IRRS, the program has been terminated with the conclusion
that the identification was stable. In this case the  potential
identified by the program was

\[
p(r)=\begin{cases}
-10.000024 & 0\leq r < 7.999994  \\
0.0 &  r \geq 7.999994
\end{cases}
\]

which is very close to the original potential

\[
q_3(r)=\begin{cases}
-10.0 & 0\leq r < 8.0\\
0.0 &  r \geq 8.0
\end{cases}
\]

On the other hand, when the phase shifts of $q_3(r)$ were corrupted by a
$10\%$ noise ($k=2.5,\ h=0.10$), the program was terminated
(according to the Stopping criterion) after 4
iterations with the Stability Index at $0.079241$. Since the Stability
Index is greater than the a priori chosen threshold of $\epsilon=0.02$
the conclusion is that the identification is unstable. A closer look
into this situation reveals that the values of the objective function
$\Phi(p_i),\ p_i\in S_{min}$ (there are 8 elements in
$S_{min}$) are between $0.0992806$ and $0.100320$.
Since we chose $\beta=1.1$ the values are within the required $10\%$ of
each other. The actual potentials for which the normalized distance is equal to
the Stability Index $0.079241$ are
\[
p_1(r)=\begin{cases}
-9.997164 & 0\leq r < 7.932678\\
-7.487082 & 7.932678 \leq r < 8.025500\\
0.0 &  r \geq 8.025500
\end{cases}
\]
and
\[
p_2(r)=\begin{cases}
-9.999565 & 0\leq r < 7.987208\\
-1.236253 & 7.987208 \leq r < 8.102628\\
0.0 &  r \geq 8.102628
\end{cases}
\]
with $\Phi(p_1)=0.0992806$ and $\Phi(p_2)=0.0997561$.
One may conclude from this example that the threshold $\epsilon=0.02$ is too tight
and can be relaxed, if the above uncertainty is acceptable.

Finally, we studied the dependency of the Stability Index from the
dimension of the admissible set $A_{adm}$, see (\ref{adm}). This
dimension is equal to $2M$ , where $M$ is the assumed number of layers
in the potential. More precisely, $M=3$, for example. means that the
search is conducted in the class of potentials having 3 or less layers.
The experiments were conducted for the identification of the original potential
$q_2(r)$ with $k=2.0$ and no noise present in the data. The results are
shown in Table \ref{table_stab5}. Since the potential $q_2$ consists of only one layer,
the smallest Stability Indices are obtained for $M=1$. They gradually
increase with $M$. Note, that the algorithm conducts the global search
using random variables, so the actual values of the indices are
different in every run. Still the results show the successful
identification (in this case) for the entire range of the a priori
chosen parameter $M$. This agrees with the theoretical consideration
according to which the Stability Index corresponding to an ill-posed
problem
in an infinite-dimensional space should be large.
Reducing the original ill-posed problem to a one in a space of much
lower dimension
regularizes the original problem.


\begin{table}
\caption{Stability Indices for $q_2(r)$ identification for different values of $M$.}

\begin{tabular}{r r r r r}
\hline

$Iteration$ & $M=1$ & $M=2$ & $M=3$ & $M=4$\\

\hline

\hline

 1 &    0.472661 & 1.068993     & 1.139720  & 1.453076 \\
 2 &    0.000000 & 0.400304     & 0.733490  & 1.453076 \\
 3 &         & 0.000426     & 0.125855  & 0.899401 \\
 4 &         &      & 0.125855  & 0.846117 \\
 5 &         &      & 0.033173  & 0.941282 \\
 6 &         &      & 0.033173  & 0.655669 \\
 7 &         &      & 0.033123  & 0.655669 \\
 8 &         &      & 0.000324  & 0.948816 \\
 9 &         &      &       & 0.025433 \\
 10 &        &      &       & 0.025433 \\
 11 &        &      &       & 0.012586 \\

\hline
\end{tabular}\label{table_stab5}
\end{table}


\section{Inverse scattering problem with fixed-energy data.}
\subsection{Problem description} In this Section we continue
 a discussion of the Inverse potential scattering with a presentation of
Ramm's method for solving inverse scattering problem with
fixed-energy data, see \cite{r474}. The method is applicable to both  exact and
noisy data. Error estimates for this method are also given.
  An inversion method using the Dirichlet-to-Neumann (DN) map
is discussed, the difficulties of its numerical implementation are pointed
out and compared with the difficulties of the implementation of the Ramm's
inversion method. See the previous Section on the potential scattering
for the problem set up.

\subsection{ Ramm's inversion method for exact data}

The results we describe in this Section are taken from \cite{r313} and
\cite{r425}.
Assume $q\in Q:=Q_a \cap L^\infty (\R^3),$ where
$Q_a:=\{ q: q(x) = \overline{q(x)}, \quad
        q(x) \in L^2(B_a), \quad q(x) = 0
        \hbox{\ if\ } |x|\geq a \},$ $B_a:=\{x: |x|\leq a\}$.
Let $A(\alpha^\prime, \alpha)$ be the corresponding scattering amplitude
at a fixed energy $k^2$, $k=1$ is taken without loss of generality. One
has:
\be\label{fix2.1}
A(\alpha^\prime, \alpha) = \sum^\infty_{\ell =0} A_\ell (\alpha)
        Y_\ell (\alpha^\prime), \quad A_\ell (\alpha) := \int_{S^2}
        A(\alpha^\prime, \alpha) \overline{Y_\ell(\alpha^\prime)}
        d \alpha^\prime,
\end{equation}
 where $S^2$ is the unit sphere in $\R^3$,
$Y_\ell (\alpha^\prime) = Y_{\ell,m} (\alpha^\prime), -\ell \leq m \leq
\ell$, are the normalized spherical harmonics,
summation over $m$ is understood in \eqref{fix2.1} and in \eqref{fix2.8} below.
Define the following algebraic variety:
\be\label{fix2.2}
  M := \{ \theta : \theta \in \C^3, \theta \cdot \theta =1\}, \quad
        \theta \cdot w := \sum^3_{j=1} \theta_j w_j.
\end{equation}
This variety is non-compact, intersects $\R^3$ over $S^2$, and, given any
$\xi \in \R^3$, there exist (many) $\theta, \theta^\prime \in M$ such that
\be\label{fix2.3}
\theta^\prime - \theta = \xi, \quad |\theta| \to \infty, \quad
        \theta, \theta^\prime \in M.
\end{equation}
In particular, if one chooses the coordinate system in which
$\xi = te_3$, $t>0$, $e_3$ is the unit vector along the $x_3$-axis, then
the vectors
\be\label{fix2.4}
        \theta^\prime = \frac{t}{2}e_3 + \zeta_2e_2 + \zeta_1e_1, \quad
        \theta = -\frac{t}{2}e_3 + \zeta_2e_2+\zeta_1 e_1, \quad
        \zeta^2_1 + \zeta^2_2 = 1-\frac{t^2}{4},
\end{equation}
satisfy \eqref{fix2.3} for any complex numbers $\zeta_1$ and $\zeta_2$ satisfying
the last
equation \eqref{fix2.4} and such that $|\zeta_1|^2+|\zeta_2|^2 \to \infty$.
There are
infinitely many such $\zeta_1, \zeta_2 \in \C$.
Consider a subset $M^\prime \subset M$ consisting of the vectors
$\theta = (\sin \vartheta \cos \varphi, \sin \vartheta \sin \varphi, \cos
\vartheta),$
where $\vartheta$ and $\varphi$ run through the whole complex plane.
Clearly
$\theta \in M$, but $M^\prime$ is a proper subset of $M$. Indeed, any
$\theta \in M$ with $\theta_3 \neq \pm 1$ is an element of $M^\prime$.
If $\theta_3 = \pm 1$, then $\cos \vartheta = \pm 1$, so
$\sin \vartheta = 0$ and one gets $\theta = (0,0, \pm 1) \in M^\prime$.
However,
there are vectors $\theta = (\theta_1, \theta_2, 1) \in M$ which do not
belong to $M^\prime$. Such vectors one obtains choosing
$\theta_1, \theta_2 \in \C$ such that $\theta^2_1 + \theta_2^2 = 0$.
There are infinitely many such vectors. The same is true for vectors
$(\theta_1, \theta_2, -1)$. Note that in \eqref{fix2.3} one can replace $M$ by
$M^\prime$ for any $ \xi \in \R^3$, $\xi \neq 2e_3$.

Let us state two estimates proved in \cite{r313}:
\be\label{fix2.5}
        \max_{\alpha \in S^2} \left| A_\ell (\alpha) \right| \leq
        c \left(\frac{a}{\ell}\right)^{\frac{1}{2}}
        \left(\frac{ae}{2\ell}\right)^{\ell +1},
\end{equation}
where $c>0$ is a constant depending on the norm $||q||_{L^2(B_a)}$, and
\be\label{fix2.6}
        \left|Y_\ell (\theta)\right| \leq \frac{1}{\sqrt{4 \pi}}
        \frac{e^{r |Im \theta|}}{|j_\ell (r)|}, \quad
        \forall r > 0, \quad \theta \in M^\prime,
\end{equation}
where
\be\label{fix2.7}
        j_\ell (r)  :=\left(\frac{\pi}{2r} \right)^{\frac{1}{2}}
        J_{\ell + \frac{1}{2}} (r) = \frac{1}{2\sqrt{2}} \frac{1}{\ell}
        \left(\frac{er}{2\ell}\right)^\ell [1 + o(1)] \hbox{\ as\ }
        \ell \to \infty,
\end{equation}
and $J_\ell(r)$ is the Bessel function regular at $r=0$.
Note that $Y_\ell (\alpha^\prime)$,
defined above, admits a natural analytic continuation from $S^2$ to
$M$
by taking $\vartheta$ and $\varphi$  to be arbitrary complex
numbers.
The resulting $\theta^\prime \in M^\prime \subset M$.

The series \eqref{fix2.1} converges absolutely and uniformly on the sets
$S^2 \times  M_c$, where $M_c$ is any compact subset of $M$.

Fix any numbers $a_1$ and $b$, such that $a<a_1<b$. Let
$||\cdot||$ denote the $L^2(a_1\leq |x|\leq b)$-norm. If $|x|\geq a,$
then the scattering solution is given analytically:
\be\label{fix2.8}
 u(x, \alpha) = e^{i \alpha \cdot x} + \sum^\infty_{\ell=0}
        A_\ell (\alpha) Y_\ell (\alpha^\prime) h_\ell (r), \quad
        r:= |x| > a, \quad \alpha^\prime := \frac{x}{r},
\end{equation}
where $A_\ell (\alpha)$ and $Y_\ell(\alpha^\prime)$
are defined above,
$$h_\ell(r) := e^{i \frac{\pi}{2} (\ell+1)} \sqrt{\frac{\pi}{2r}}
H^{(1)}_{\ell + \frac{1}{2}} (r),$$
$H^{(1)}_\ell (r)$ is the Hankel function, and the normalizing factor is
chosen so that
$        h_\ell(r) = \frac{e^{ir}}{r} [1 + o(1)] \hbox{\ as\ } r \to
\infty.$
Define
\be\label{fix2.9}
 \rho(x) := \rho(x;\nu) := e^{-i \theta \cdot x} \int_{S^2}
        u(x, \alpha) \nu(\alpha, \theta)d \alpha -1, \quad \nu \in
L^2(S^2).
\end{equation}
Consider the minimization problem
\be\label{fix2.10}
        \Vert \rho \Vert = \inf := d(\theta),
\end{equation}
where the infimum is taken over all $\nu \in L^2(S^2)$, and \eqref{fix2.3} holds.

It is proved in \cite{r313}  that
\be\label{fix2.11}
        d(\theta) \leq c |\theta|^{-1} \hbox{\ if\ } \theta \in M, \quad
        |\theta| \gg 1.
\end{equation}
The symbol $|\theta| \gg 1$ means that $|\theta|$ is sufficiently large.
The
constant $c>0$ in \eqref{fix2.11} depends on the norm $\Vert q \Vert_{L^2(B_a)}$
but
not on the potential $q(x)$ itself.

{\it An algorithm for computing a function
$\nu(\alpha, \theta)$, which can be used for inversion of the
exact, fixed-energy, three-dimensional
scattering data, is as follows:}

a) Find an approximate solution to \eqref{fix2.10} in the sense
\be\label{fix2.12}
        \Vert \rho (x, \nu) \Vert < 2 d(\theta),
\end{equation}
where in place of the factor 2 in \eqref{fix2.12} one could put any fixed constant
greater than 1.

b) Any such $\nu (\alpha, \theta)$ generates an estimate of
$\widetilde {q}(\xi)$ with the error $O\left(\frac{1}{|\theta|}\right)$,
$|\theta | \to \infty$.
This estimate is calculated by the formula
\be\label{fix2.13}
        \widehat{q} := -4\pi \int_{S^2} A(\theta^\prime, \alpha)
        \nu (\alpha, \theta) d \alpha,
\end{equation}
where $\nu(\alpha, \theta) \in L^2(S^2)$ is any function satisfying
\eqref{fix2.12}.

Our basic result is:

\begin{theorem}\label{fixthm1}  Let \eqref{fix2.3} and \eqref{fix2.12} hold. Then
\be\label{fix2.14}
    \sup_{\xi \in \R^3}   \left|\widehat{q} - \widetilde{q}(\xi)
\right| \leq
        \frac{c}{|\theta|}, \quad
        |\theta| \to \infty, \quad
\end{equation}
The constant $c>0$ in \eqref{fix2.14} depends
on a norm of $q$, but not on a particular $q$.
\end{theorem}

The norm of $q$ in the above Theorem can be any norm such that the set
$\{q: ||q||\leq const\}$ is a compact set in $L^\infty(B_a)$.

In \cite{r313} and \cite{r425} an inversion algorithm is formulated also for
noisy data, and the error estimate for this algorithm is obtained.
 Let us describe these results.

Assume that the scattering data are given with some error: a function
$A_\delta (\alpha^\prime, \alpha)$ is given such that
\be\label{fix2.15}
        \sup_{\alpha^\prime, \alpha \in S^2}
        \left|A(\alpha^\prime, \alpha) -
        A_\delta(\alpha^\prime, \alpha) \right| \leq \delta.
\end{equation}

We emphasize that $A_\delta (\alpha^\prime, \alpha)$ is not necessarily a
scattering amplitude corresponding to some potential, it is an arbitrary
 function in $L^\infty (S^2 \times S^2)$ satisfying \eqref{fix2.15}. It is assumed
that
the unknown function $A(\alpha^\prime, \alpha)$ is the scattering
amplitude
corresponding to a $q \in Q$.

The problem is: {\it Find an algorithm for calculating
$\widehat{q_\delta}$
such that}
\be\label{fix2.16}
        \sup_{\xi \in \R^3} \left|\widehat{q_\delta} - \widetilde{q}
(\xi)
        \right| \leq \eta(\delta), \quad \eta(\delta) \to 0
        \hbox{\ as\ } \delta \to 0,
\end{equation}
{\it  and estimate
the rate at which $\eta(\delta)$ tends to zero.}

An algorithm for inversion of noisy data will now be described.

Let
\be\label{fix2.17}
        N(\delta) := \left[ \frac{|\ln \delta|}{\ln|\ln \delta|} \right],
\end{equation}
where $[x]$ is the integer nearest to $x>0$,
\be\label{fix2.18}
        \widehat{A}_\delta (\theta^\prime, \alpha) :=
        \sum^{N(\delta)}_{\ell =0} A_{\delta \ell} (\alpha) Y_\ell
        (\theta^\prime), \quad
        A_{\delta \ell} (\alpha) := \int_{S^2}
        A_\delta (\alpha^\prime, \alpha) \overline{Y_\ell(\alpha^\prime)}
        d \alpha^\prime,
\end{equation}
\be\label{fix2.19}
        u_\delta (x, \alpha) := e^{i \alpha \cdot x} +
        \sum^{N(\delta)}_{\ell =0} A_{\delta \ell}
        (\alpha) Y_\ell (\alpha^\prime) h_\ell (r),
\end{equation}
\be\label{fix2.20}
        \rho_\delta (x; \nu) := e^{-i \theta \cdot x} \int_{S^2}
        u_\delta (x, \alpha) \nu (\alpha) d \alpha -1, \quad
        \theta \in M,
\end{equation}
\be\label{fix2.21}
        \mu(\delta) := e^{-\gamma N(\delta)}, \quad
        \quad \gamma=\ln \frac {a_1}{a} >0,
\end{equation}
\be\label{fix2.22}
        a(\nu) := \Vert \nu \Vert_{L^2(S^2)}, \quad
        \kappa := |Im \theta|.
\end{equation}
Consider the variational problem with constraints:
\be\label{fix2.23}
        |\theta| = \sup := \vartheta(\delta),
\end{equation}
\be\label{fix2.24}
       |\theta| \left[ \Vert \rho_\delta (\nu) \Vert +
        a(\nu) e^{\kappa b} \mu(\delta) \right] \leq c,
        \quad \theta \in M, \quad
        |\theta| = \sup := \vartheta(\delta),
\end{equation}
the norm is defined above \eqref{fix2.8}, and it is assumed that \eqref{fix2.3} holds,
where $\xi \in \R^3$ is an arbitrary fixed vector, $c>0$ is a sufficiently
large constant, and the supremum is taken over $\theta \in M$ and
$\nu \in L^2(S^2)$ under the constraint \eqref{fix2.24}. By $c$ we denote various
positive constants.

Given $\xi \in \R^3$ one can always find $\theta$ and $\theta^\prime$
such that \eqref{fix2.3} holds.
  We prove that $\vartheta(\delta) \to \infty$, more precisely:
\be\label{fix2.25}
        \vartheta(\delta) \geq c \frac{|\ln \delta|}{(\ln|\ln
\delta|)^2},
        \quad \delta \to 0.
\end{equation}

Let the pair
$\theta(\delta)$ and $\nu_\delta (\alpha, \theta)$ be any approximate
solution to
problem \eqref{fix2.23}-\eqref{fix2.24} in the sense that
\be\label{fix2.26}
        |\theta(\delta)| \geq \frac{\vartheta(\delta)}{2}.
\end{equation}
Calculate
\be\label{fix2.27}
        \widehat{q}_\delta := -4 \pi \int_{S^2} \widehat{A}_\delta
        (\theta^\prime, \alpha) \nu_\delta (\alpha, \theta) d \alpha.
\end{equation}

\begin{theorem}\label{fixthm2}
      If \eqref{fix2.3} and \eqref{fix2.26} hold, then
\be\label{fix2.28}
        \sup_{\xi \in \R^3} \left| \widehat{q}_\delta -
        \widetilde{q}(\xi) \right| \leq c
        \frac{(\ln |\ln \delta|)^2}{|\ln \delta|} \hbox{\ as\ }
        \delta \to 0,
\end{equation}
where $c>0$ is a constant depending on a norm of $q$.
\end{theorem}
 In \cite{r313} estimates \eqref{fix2.14} and \eqref{fix2.28} were formulated with the supremum
taken over an arbitrary large but fixed ball of radius $\xi_0$. Here these
estimates are improved: $\xi_0=\infty$. The key point is: the constant
$c>0$ in the estimate \eqref{fix2.11} does not depend on 
$\theta$.

{\bf Remark.} In \cite{R8} (see also \cite{R9} and \cite{r425}) an analysis
of the approach to ISP, based on the recovery of the DN
(Dirichle-to-Neumann) map from the fixed-energy scattering data, is given.
This approach is discussed below.

The basic numerical difficulty of the approach described in Theorems
\ref{fixthm1}
and \ref{fixthm2} comes from solving problems \eqref{fix2.10} for exact data, and
problem \eqref{fix2.23}-\eqref{fix2.24}  for noisy data. Solving \eqref{fix2.10} amounts to
finding a global minimizer of a quadratic form of the variables $c_\ell$,
if one
takes $\nu$ in \eqref{fix2.9} as a linear combination of the spherical harmonics:
$\nu=\sum_{\ell =0}^L c_\ell  Y_\ell (\alpha)$. If one uses
the necessary condition for a minimizer of a quadratic form, that is, a
linear system, then the matrix of this system is ill-conditioned
for large $L$. This causes the main difficulty in the numerical
solution of \eqref{fix2.10}. On the other hand, there are methods for global
minimization of the quadratic functionals, based on the gradient descent,
which may be more efficient than using the above necessary condition.

\subsection{Discussion of the inversion method which uses the DN map}

In \cite{R8} the following inversion method is discussed:

\be\label{fix3.1}
\tilde q (\xi)=\lim_{|\theta|\to \infty} \int_S\exp
(-i\theta'\cdot s)(\Lambda-\Lambda_0)\psi ds ,
\end{equation}
where \eqref{fix2.3} is assumed, $\Lambda$ is the Dirichlet-to-Neumann (DN) map,
$\psi$ is
found from the equation:
\be\label{fix3.2}
\psi (s)=\psi_0(s)-\int_S G(s-t)B\psi dt,\quad B:=\Lambda-\Lambda_0,
\quad \psi_0(s):=e^{i\theta \cdot s},
\end{equation}
and $G$ is defined by the formula:
\be\label{fix3.3}
G(x)=\exp (i\theta\cdot x)\frac 1{(2\pi)^3}\int_{\mathbb R^3}
\frac{\exp (i\xi\cdot x)d\xi}{\xi^2+2\xi\cdot\theta}.
\end{equation}
The DN map is constructed from the
fixed-energy scattering data $A(\alpha', \alpha)$
by the method of \cite{R8} (see also \cite{r313}).

Namely, given $A(\alpha', \alpha)$ for all $\alpha', \alpha \in S^2$,
one finds $\Lambda$ using the following steps.

Let $f\in H^{3/2}(S)$ be given, $S$ is a sphere of radius $a$ centered
at the origin, $f_\ell$ are its Fourier coefficients
in the basis of the spherical harmonics,
\be\label{fix3.4}
w=\sum^\infty_{l=0}f_lY_l(x^0)\frac{h_l(r)}{h_l(a)},\quad r\geq a,
\quad x^0:=\frac xr, \quad r:=|x|.
\end{equation}
Let
\be\label{fix3.5}
w=\int_S g(x,s)\sigma (s)ds ,
\end{equation}
where $\sigma$ is some function, which we find below, and $g$ is the
Green function (resolvent kernel) of the Schroedinger operator,
satisfying the radiation condition at infinity.
Then
\be\label{fix3.6}
w^+_N=w^-_N+\sigma ,
\end{equation}
where $N$ is the outer normal to $S$, so $N$ is directed along the
radius-vector.
We require $w=f$ on $S$. Then $w$ is given by \eqref{fix3.4} in the exterior of
$S$, and
\be\label{fix3.7}
w_N^-=\sum^\infty_{l=0}f_lY_l(x^0)\frac{h^\prime_l(a)}{h_l(a)}.
\end{equation}
By formulas \eqref{fix3.6} and \eqref{fix3.7}, finding $\Lambda$ is equivalent to finding
$\sigma$. By \eqref{fix3.5}, asymptotics of $w$ as $r:=|x|\to \infty$,
$x/|x|:=x^0,$ is (cf \cite{r313}, p.67):
\be\label{fix3.8}
w=\frac {e^{ir}}{r} \frac {u(y,-x^0)}{4\pi} +o(\frac 1 r),
\end{equation}
where $u$ is the scattering solution,
\be\label{fix3.9}
u(y,-x^0)=e^{-ix^0 \cdot y}+\sum_{\ell=0}^\infty A_\ell
(-x^0)Y_\ell(y^0)h_\ell(|y|).
\end{equation}
>From \eqref{fix3.4}, \eqref{fix3.8} and \eqref{fix3.9} one gets an equation for finding $\sigma$
(\cite{R8}, eq. (23), see also \cite{r313}, p. 199):
\be\label{fix3.10}
\frac{f_l}{h_l(a)}=\frac 1{4\pi}\int_S ds\sigma (s)\left (u(s,-\beta),
Y_l(\beta)\right )_{L^2(S^2)} ,
\end{equation}
which can be written as a linear system:
\be\label{fix3.11}
\frac{4\pi f_l}{h_l(a)}=a^2(-1)^l\sum^\infty_{l^\prime
=0}\sigma_{l^\prime}
[4\pi i^lj_l(a)\delta_{ll^\prime}+A_{l^\prime l}h_{l^\prime}(a)] ,
\end{equation}
for the Fourier coefficients $\sigma_\ell$ of $\sigma$. The coefficients
$$A_{l^\prime l}:=((A(\alpha', \alpha), Y_\ell(\alpha'))_{L^2(S^2)}, Y_\ell
(\alpha))_{L^2(S^2)}
$$
are the Fourier coefficients of the scattering
amplitude. Problems \eqref{fix3.10} and \eqref{fix3.11} are very ill-posed
(see \cite{R8} for details).

{\it This approach faces many difficulties:}

1) The construction of the DN map from the scattering data is a very
ill-posed problem,

2) The construction of the potential
from the DN map is a very difficult problem numerically,  because
one has to solve a Fredholm-type integral equation
( equation \eqref{fix3.2} ) whose kernel contains $G$, defined in \eqref{fix3.3}. This $G$
is a
tempered distribution, and it is
very difficult to compute it,

3) One has to calculate a limit
of an integral whose integrand grows exponentially to infinity
if a factor in the integrand is not known exactly. The solution of
equation \eqref{fix3.2}  is one of the factors in the integrand.
It cannot be known exactly in practice because it cannot be
calculated with arbitrary accuracy even if the scattering data are known
exactly. Therefore the limit in formula \eqref{fix3.1}  cannot be
calculated accurately.

No error estimates are obtained for this approach.

In contrast, in Ramm's method, there is no need to
compute $G$, to solve equation \eqref{fix3.2}, to calculate the DN map from
the scattering data, and to compute the limit \eqref{fix3.1}.
The basic difficulty in Ramm's inversion method for exact data is to
minimize the quadratic form \eqref{fix2.10}, and for noisy data to solve optimization
problem \eqref{fix2.23}-\eqref{fix2.24}. The error estimates are obtained for the Ramm's
method.

\section{Obstacle scattering by the Modified Rayleigh Conjecture  (MRC) method.}

\subsection{Problem description}
In this section we present a novel numerical method for Direct
Obstacle
 Scattering  Problems based on the Modified Rayleigh Conjecture
(MRC). The basic theoretical foundation of the method was developed in
\cite{r430}. The MRC has the appeal of an easy implementation for obstacles
of complicated geometry, e.g. having edges and corners. A special version of the MRC method
was used in \cite{r475} to compute the scattered field for 3D obstacles.
In our numerical
experiments the method has shown itself to be a competitive alternative to the BIEM (boundary
integral equations method), see \cite{r437}. Also, unlike the BIEM, one
can apply the algorithm to different obstacles with very little
additional effort.

We formulate the obstacle scattering problem in a 3D setting with the
 Dirichlet boundary condition, but the
discussed method can also be used for the Neumann and Robin boundary
conditions.

Consider a bounded domain $D \subset \R^3$, with a boundary $S$ which is
assumed to be Lipschitz continuous. Denote the
exterior domain by $D^\prime = \R^3 \backslash D$. Let
 $\alpha, \alpha^\prime\in S^2$ be unit vectors,
and $ S^2$ be the unit sphere in $\R^3$.

The acoustic wave scattering problem by a soft obstacle $D$
consists in
finding the (unique) solution to the problem \eqref{s1_e1}-\eqref{s1_e2}:
\be\label{s1_e1}
\left(\nabla^2 + k^2 \right) u=0 \hbox{\ in\ } D^\prime, \quad
  u = 0 \hbox{\ on\ } S,
\end{equation}

\be\label{s1_e2}
u=u_0 + A(\alpha^\prime, \alpha) \frac{e^{ikr}}{r}
+ o
  \left(\frac{1}{r} \right), \quad r:=|x| \to \infty, \quad
  \alpha^\prime := \frac{x}{r}.
\end{equation}
Here $u_0:=e^{ik \alpha \cdot x}$ is the incident field,
$v:=u-u_0$ is the scattered field,
$A(\alpha^\prime, \alpha)$ is called the scattering amplitude, its
k-dependence is not shown, $k>0$ is the wavenumber. Denote

\be\label{s1_e3}
A_\ell (\alpha) := \int_{S^2} A(\alpha^\prime, \alpha)
  \overline{Y_\ell (\alpha^\prime)} d\alpha^\prime,
\end{equation}
where $Y_\ell (\alpha)$ are the orthonormal spherical harmonics,
$Y_\ell = Y_{\ell m}, -\ell \leq m \leq \ell$. Let $h_\ell (r)$ be the spherical
Hankel functions, normalized so that
$h_\ell (r) \sim \frac{e^{ikr}}{r}$ as $r \to +\infty$.

Informally, the Random Multi-point MRC algorithm can be described as follows.

Fix a $J>0$. Let $x_j, j=1,2,...,J$ be a batch of points randomly chosen  inside the
obstacle $D$. For $x\in D^\prime$, let
\begin{equation}\label{s1_hy}
\alpha^\prime =
  \frac{x-x_j}{|x-x_j|},\quad \psi_\ell(x,x_j)= Y_\ell
(\alpha^\prime) h_\ell (k|x-x_j|).  \end{equation} Let
$g(x)=u_0(x),\; x\in S$, and minimize the discrepancy
\be\label{s1_phi} \Phi(\bc)=\|g(x)+\sum_{j=1}^J\sum_{\ell=0}^L
c_{\ell,j}\psi_\ell(x,x_j)\|_{L^2(S)}, \end{equation}
 over
$\bc\in \C^N$, where $\bc=\{c_{\ell,j}\}$. That is, the total
field $u=g(x)+v$ is desired to be as close to zero as
possible at the boundary $S$, to satisfy the required
condition for the soft scattering. If the resulting residual
$r^{min}=\min\Phi$ is smaller than the prescribed tolerance
$\ep$, than the procedure is finished, and the sought
scattered field is
$$v_{\epsilon}(x)=
\sum_{j=1}^J\sum_{\ell=0}^{L}c_{\ell,j}\psi_\ell(x,x_j),\quad
x\in D',$$
(see Lemma \ref{s1_lm1} below).

If, on the other hand, the residual $r^{min}>\ep$, then we
continue by trying to improve on the already obtained fit in
(\ref{s1_phi}). Adjust the field on the boundary by letting
$g(x):=g(x)+v_{\epsilon}(x),\; x\in S$. Create another batch
of $J$ points randomly chosen in the interior of $D$, and
minimize (\ref{s1_phi}) with this new $g(x)$. Continue with
the iterations until the required tolerance $\ep$ on the
boundary $S$ is attained, at the same time keeping the track
of the changing field $v_{\epsilon}$.

Note, that the minimization in (\ref{s1_phi}) is always done
over the same number of points $J$. However, the points
$x_j$ are sought to be different in each iteration
 to assure that the minimal values of $\Phi$ are decreasing
in consequent iterations. Thus, computationally, the size of
the minimization problem remains the same. This is the new
feature of the Random multi-point MRC method, which allows
it to solve scattering problems untreatable by previously
developed MRC methods, see \cite{r437}.

Here is the precise description of the algorithm.

{\bf Random Multi-point MRC.}

For $x_j\in D$, and $\ell\geq 0$ functions $\psi_\ell(x,x_j)$
are defined as in (\ref{s1_hy}).

\begin{enumerate} \item\label{init4} {\bf Initialization.}
Fix $\ep>0, \; L\geq 0,\;J>0,\; N_{max}>0$. Let $n=0$,
$v_{\ep}=0$ and $g(x)=u_0(x), \;x\in S$.

\item\label{iter4} {\bf Iteration.}

\begin{enumerate}
 \item Let $n:=n+1$. Randomly choose $J$ points $x_j\in D,\; j=1,2,\dots, J$.

 \item
Minimize
\[
\Phi(\bc)=\|g(x)+\sum_{j=1}^J\sum_{\ell=0}^L c_{\ell,j}\psi_\ell(x,x_j)\|_{L^2(S)}
\]
over $\bc\in \C^N$, where $\bc=\{c_{\ell,j}\}$.

Let the minimal value of $\Phi$ be $r^{min}$.

 \item Let
\[
v_{\epsilon}(x):=v_{\epsilon}(x)+\sum_{j=1}^J\sum_{\ell=0}^{L}
c_{\ell,j}\psi_\ell(x,x_j),\quad x\in D'.
\]

\end{enumerate}

\item {\bf Stopping criterion.}
\begin{enumerate}
 \item If $r^{min}\leq\epsilon$, then  stop.

 \item If $r^{min}>\epsilon$, and $n\not=N_{max}$,
let
\[
g(x):=g(x)+\sum_{j=1}^J\sum_{\ell=0}^{L}c_{\ell,j}\psi_\ell(x,x_j),\quad x\in S
\]
and repeat the iterative
step (\ref{iter4}).
\item If $r^{min}>\ep$, and $n=N_{max}$, then the procedure failed.

\end{enumerate}
\end{enumerate}

\subsection{Direct scattering problems and the Rayleigh conjecture.}

Let a ball
$B_R := \{x : |x| \leq R\}$ contain the obstacle $D$.
In the region $r>R$ the solution to (\ref{s1_e1})-(\ref{s1_e2}) is:
\be\label{ss2_e1}
u(x, \alpha) = e^{ik\alpha \cdot x} +
\sum^\infty_{\ell =0} A_\ell (\alpha)\psi_\ell, \quad
\psi_\ell:= Y_\ell (\alpha^\prime) h_\ell (kr), \quad r > R,\quad
\alpha^\prime =
  \frac{x}{r},
\end{equation}
where the sum includes the summation with respect to $m$, $-\ell \leq m \leq \ell$,
and $A_\ell (\alpha)$ are defined in (\ref{s1_e3}).

{\it The Rayleigh conjecture (RC) is:  the series (\ref{ss2_e1}) converges up to
the
boundary $S$} (originally RC dealt with periodic structures, gratings).
This conjecture is false for many obstacles, but is true for some
(\cite{baran,millar,r190}). For example, if $n=2$
and
$D$ is an ellipse, then the series analogous to (\ref{ss2_e1}) converges in the
region
$r >a$, where $2a$ is the distance between the foci of the ellipse \cite{baran}.
In the engineering literature there are numerical algorithms, based on the
Rayleigh conjecture.
Our aim is to give a formulation of a {\it Modified Rayleigh Conjecture}
(MRC)
which holds for any Lipschitz obstacle and can be used in numerical
solution of the direct and
inverse scattering problems (see \cite{r430}). We discuss the Dirichlet
condition but
similar argument is applicable to the Neumann boundary
condition, corresponding to acoustically hard obstacles.

Fix $\epsilon >0$, an arbitrary small number.

\begin{lemma}\label{s1_lm1}
There exist $L=L(\epsilon)$ and
$c_\ell=c_\ell(\epsilon)$
such that
\be\label{ss2_e6}
||u_0+\sum_{\ell=0}^{L(\epsilon)}c_\ell(\epsilon)\psi_\ell||_{L^2(S)} \leq
\epsilon.
\end{equation}

If (\ref{ss2_e6}) and the boundary condition (\ref{s1_e1}) hold, then
\be\label{ss2_e7}
||v_{\epsilon}-v||_{L^2(S)}\leq \epsilon,  \quad
v_{\epsilon}:=\sum_{\ell=0}^{L(\epsilon)}c_\ell(\epsilon)\psi_\ell.
\end{equation}
\end{lemma}

\begin{lemma}\label{s1_lm2}
 If (\ref{ss2_e7}) holds then
\be\label{ss2_e8}
|||v_{\epsilon}-v|||=O(\epsilon)\,, \quad \epsilon \to 0,
\quad
\end{equation}
where $|||\cdot|||:=
||\cdot||_{H_{loc}^m(D')}+||\cdot||_{L^2(D';
(1+|x|)^{-\gamma})}$, $\gamma >1$, $m>0$ is an arbitrary integer,
$H^m$ is the Sobolev space, and $v_{\epsilon}, v$ in (\ref{ss2_e8}) are
functions defined in $D'$.

In particular, (\ref{ss2_e8}) implies
\be\label{ss2_e9}
||v_{\epsilon}-v||_{L^2(S_R)}=O(\epsilon)\,, \quad  \epsilon \to 0,
\end{equation}
where $S_R$ is the sphere centered at the origin with radius $R$.
\end{lemma}

\begin{lemma}\label{s1_lm3}
 One has:
\be\label{ss2_e10}
 c_\ell(\epsilon) \to A_\ell(\alpha) \,,\quad \forall \ell, \quad \epsilon \to 0.
\end{equation}
\end{lemma}

The Modified Rayleigh Conjecture (MRC) is formulated as a theorem, which
follows from the above three lemmas:

\begin{theorem}\label{s1_thm1}
 For an arbitrary small $\epsilon>0$ there
exist
$L(\epsilon)$ and $c_\ell(\epsilon),\,\, 0\leq \ell \leq L(\epsilon)$,
such that (\ref{ss2_e6}), (\ref{ss2_e8}) and (\ref{ss2_e10}) hold.
\end{theorem}

See \cite{r430} for a proof of the above statements.

The difference between RC and MRC is: (\ref{ss2_e7}) does
not hold if one replaces $v_\epsilon$ by $\sum_{\ell=0}^L
A_\ell(\alpha)\psi_\ell$, and lets $L\to \infty$ (instead of
letting $\epsilon \to 0$). Indeed, the series
$\sum_{\ell=0}^\infty A_\ell(\alpha)\psi_\ell$ diverges at some
points of the boundary for many obstacles. Note also that
the coefficients in (\ref{ss2_e7}) depend on $\epsilon$, so
(\ref{ss2_e7}) is {\it not} a partial sum of a series.

For the Neumann boundary condition one minimizes
$$ \left\|\frac {\partial [u_0+\sum_{\ell=0}^{L}c_\ell\psi_\ell]}{\partial
N}\right\|_{L^2(S)}$$
with respect to $c_\ell$. Analogs of Lemmas \ref{s1_lm1}-\ref{s1_lm3}
are valid and their
proofs are essentially the same.

See \cite{r461} for an extension of these results to scattering by periodic
structures.

\subsection{Numerical Experiments.}

In this section we desribe numerical results obtained by the Random
Multi-point MRC method for 2D and 3D obstacles. We also compare the 2D results
 to the ones obtained by our earlier method introduced in
\cite{r437}. The method that we used previously can be
described as a Multi-point MRC. Its difference from the
Random Multi-point MRC method is twofold: It is just the
first iteration of the Random method, and the interior
points $x_j,\; j=1,2,...,J$ were chosen deterministically,
by an {\it ad hoc} method
 according to the geometry of the obstacle $D$. The number
of points $J$ was limited by the size of the resulting
numerical minimization problem, so the accuracy of the
scattering solution (i.e. the residual $r^{min}$) could not
be made small for many obstacles. The method was not capable
of treating 3D obstacles. These limitations were removed by
using the Random Multi-point MRC method. As we mentioned
previously, \cite{r437} contains a favorable comparison of
the Multi-point MRC method with the BIEM, inspite in spite
of the fact that the numerical implementation of the MRC
method in \cite{r437} is considerably less efficient than the one
presented in this paper.

A numerical implementation of the Random Multi-point MRC
method follows the same outline as for the Multi-point MRC,
which was described in \cite{r437}. Of course, in a 2D case,
instead of (\ref{s1_hy}) one has

\[
\psi_l(x,x_j)=H_l^{(1)}(k\abs{x-x_j})e^{il\theta_j},
\]
where $(x-x_j)/\abs{x-x_j}=e^{i\theta_j}$.

For a  numerical implementation
choose
$M$ nodes $\{t_m\}$ on the surface $S$ of the obstacle $D$. After the interior points
$x_j,\; j=1,2,...,J$ are chosen, form
 $N$ vectors
\[
{\bf a}^{(n)}=\{\psi_l(t_m,x_j)\}_{m=1}^M,
\]
$n=1,2,\dots,N$ of length $M$. Note that $N=(2L+1)J$
for a 2D case, and $N=(L+1)^2J$ for a 3D case.
It is convenient to normalize the norm in $\R^M$ by
\[
\|{\bf b}\|^2=\frac 1M \sum_{m=1}^M|b_m|^2,\quad {\bf
b}=(b_1,b_2,...,b_M).
\]
Then $\|u_0\|=1$.

Now let ${\bf b}=\{g(t_m)\}_{m=1}^M$, in
the Random Multi-point MRC (see section 1), and minimize
\be\label{s3_minm}
\Phi({\bf c})=\|{\bf b}+A{\bf c}\|,
\end{equation}
for ${\bf c}\in \C^N$, where $A$ is the matrix containing vectors ${\bf a}^{(n)},\;
n=1,2,\dots,N$ as its columns.

We used the Singular Value Decomposition (SVD) method (see e.g.
\cite{numrec}) to minimize (\ref{s3_minm}).  Small singular values $s_n<w_{min}$ of the
matrix $A$ are used to identify and delete linearly dependent or almost
linearly dependent combinations of vectors ${\bf a}^{(n)}$. This spectral
cut-off makes the minimization process stable, see the details in
\cite{r437}.

{\it Let $r^{min}$ be the residual, i.e. the minimal value
of
$\Phi({\bf c})$ attained after $N_{max}$ iterations of the
Random Multi-point MRC method (or when it is stopped). For a
comparison, let $r^{min}_{old}$ be the residual obtained in
\cite{r437} by an earlier method.}

We conducted 2D numerical experiments for four obstacles: two ellipses of
different eccentricity, a kite, and a triangle.
The M=720 nodes $t_m$ were uniformly distributed on the interval
$[0,2\pi]$, used to parametrize the boundary $S$.
Each case was tested for
wave numbers $k=1.0$ and $ k=5.0$. Each obstacle was subjected to incident
waves corresponding to $\a=(1.0,0.0)$ and $\a=(0.0,1.0)$.

The results for the Random Multi-point MRC with $J=1$ are shown in Table \ref{mrc_table1}, in the
last column $r^{min}$. In every experiment the target
residual $\ep=0.0001$ was obtained in under 6000 iterations,
 in about 2 minutes run time on a 2.8 MHz PC.

In \cite{r437}, we conducted numerical experiments for the same four 2D
obstacles by a Multi-point MRC, as described in the beginning of this
section. The interior points $x_j$ were chosen differently in
each experiment. Their choice is indicated in the description of each 2D
experiment.
The column $J$ shows the number of these interior
points. Values $L=5$ and $M=720$ were used in all the
experiments. These results are
shown in Table \ref{mrc_table1}, column  $r^{min}_{old}$.

 Thus, the Random Multi-point MRC method achieved a significant
improvement over the earlier Multi-point MRC.

\begin{table}
\caption{Normalized residuals attained in the
numerical experiments for 2D obstacles,
 $\|{\bf u_0}\|=1$.}

\begin{tabular}{c  r  r  c  r   r}

\hline
Experiment & $J$ & $k$ & $\a$ & $r^{min}_{old}$ & $r^{min}$ \\

\hline
I & 4 & 1.0 & $(1.0,0.0)$ & 0.000201  & 0.0001\\
  & 4 & 1.0 & $(0.0,1.0)$ & 0.000357  & 0.0001\\
  & 4 & 5.0 & $(1.0,0.0)$ & 0.001309  & 0.0001\\
  & 4 & 5.0 & $(0.0,1.0)$ & 0.007228  & 0.0001\\
\hline
II & 16 & 1.0 & $(1.0,0.0)$ & 0.003555  & 0.0001\\
  & 16 & 1.0 & $(0.0,1.0)$ & 0.002169  & 0.0001\\
  & 16 & 5.0 & $(1.0,0.0)$ & 0.009673  & 0.0001\\
  & 16 & 5.0 & $(0.0,1.0)$ & 0.007291  & 0.0001\\
\hline
III & 16 & 1.0 & $(1.0,0.0)$ & 0.008281  & 0.0001\\
    & 16 & 1.0 & $(0.0,1.0)$ & 0.007523  & 0.0001\\
    & 16 & 5.0 & $(1.0,0.0)$ & 0.021571  & 0.0001\\
    & 16 & 5.0 & $(0.0,1.0)$ & 0.024360  & 0.0001\\
\hline
IV & 32 & 1.0 & $(1.0,0.0)$ & 0.006610  & 0.0001\\
   & 32 & 1.0 & $(0.0,1.0)$ & 0.006785  & 0.0001\\
   & 32 & 5.0 & $(1.0,0.0)$ & 0.034027  & 0.0001\\
   & 32 & 5.0 & $(0.0,1.0)$ & 0.040129  & 0.0001\\
\hline
\end{tabular}\label{mrc_table1}

\end{table}

{\bf Experiment 2D-I.} The boundary $S$ is an ellipse described by
\be
{\bf r}(t)=(2.0\cos t,\ \sin t),\quad 0\leq t<2\pi\,.
\end{equation}
The Multi-point MRC used $J=4$ interior points
$x_j=0.7{\bf r}(\frac{\pi(j-1)}2),\; j=1,\dots,4$.
Run time was 2 seconds.

{\bf Experiment 2D-II.} The kite-shaped boundary $S$ (see \cite{coltonkress}, Section 3.5) is
described by
\be
{\bf r}(t)=(-0.65+\cos t+0.65\cos 2t,\ 1.5\sin t),\quad 0\leq t<2\pi\,.
\end{equation}
The Multi-point MRC used $J=16$ interior points $x_j=0.9{\bf r}(\frac{\pi(j-1)}8),\; j=1,\dots,16$.
Run time  was 33 seconds.


{\bf Experiment 2D-III.} The boundary $S$ is the triangle
with vertices $(-1.0,0.0)$ and $(1.0,\pm 1.0)$.
The Multi-point MRC used the interior points $x_j=0.9{\bf r}(\frac{\pi(j-1)}8)$,
$ j=1,\dots,16$.
Run time  was about 30 seconds.

{\bf Experiment 2D-IV.} The boundary $S$ is an ellipse described by
\be
{\bf r}(t)=(0.1\cos t,\ \sin t),\quad 0\leq t<2\pi\,.
\end{equation}
The Multi-point MRC used $J=32$ interior points
$x_j=0.95{\bf r}(\frac{\pi(j-1)}{16}),\; j=1,\dots,32$.
Run time was about 140 seconds.

The 3D numerical experiments were conducted for 3 obstacles: a sphere, a cube, and an ellipsoid.
We used the Random Multi-point MRC with $L=0,\; w_{min}=10^{-12}$, and $J=80$.
 The number $M$ of the points on the boundary $S$ is
indicated in the description of the obstacles. The scattered
field for each obstacle was computed for two incoming
directions $\a_i=(\theta,\phi),\;i=1,2$, where $\phi$ was
the polar angle. The first unit vector $\a_1$ is denoted by
(1) in Table \ref{mrc_table2}, $\a_1=(0.0,\pi/2)$. The second one is
denoted by (2), $\a_2=(\pi/2,\pi/4)$. A typical number of
iterations $N_{iter}$ and the run time on a 2.8 MHz PC are
also shown in Table \ref{mrc_table2}. For example, in experiment I with
$k=5.0$ it took about 700 iterations of the Random
Multi-point MRC method to achieve the target residual
$r^{min}=0.001$ in 7 minutes.

{\bf Experiment 3D-I.} The boundary $S$ is the sphere of radius $1$,
with $M=450$.

{\bf Experiment 3D-II.} The boundary $S$ is the surface of the cube
$[-1,1]^3$ with $M=1350$.

{\bf Experiment 3D-III.} The boundary $S$ is the surface of the
ellipsoid $x^2/16+y^2+z^2=1$ with $M=450$.

\begin{table}
\caption{Normalized residuals attained in the numerical experiments for 3D obstacles,
$\|{\bf u_0}\|=1$.}

\begin{tabular}{c  r  r  l  r   r}

\hline
Experiment & $k$ & $\a_i$ &  $r^{min}$ & $N_{iter}$ & run time \\

\hline
I & 1.0 &  & $0.0002$ & 1   & 1 sec\\
  & 5.0 &  & $0.001$ & 700  & 7 min\\

\hline
II & 1.0 & (1) & $0.001$ & 800   & 16 min\\
   & 1.0 & (2) & $0.001$ & 200  & 4 min\\
   & 5.0 & (1) & $0.0035$ & 2000   & 40 min\\
   & 5.0 & (2) & $0.002$ & 2000  & 40 min\\

\hline
III & 1.0 & (1) & $0.001$ & 3600   & 37 min\\
    & 1.0 & (2) & $0.001$ & 3000   & 31 min\\
    & 5.0 & (1) & $0.0026$ & 5000   & 53 min\\
    & 5.0 & (2) & $0.001$ & 5000   & 53 min\\
\hline
\end{tabular}\label{mrc_table2}

\end{table}
In the last experiment the run time could be reduced by taking a smaller value
for $J$. For example, the choice of $J=8$ reduced the running time to
about 6-10 minutes.

Numerical experiments show that the minimization results depend on the
choice of such parameters as $J,\; w_{min}$, and $L$. They also depend on
the choice of the interior points $x_j$. It is possible that further
versions of the MRC could be made more efficient by finding a more
efficient rule for their placement. Numerical experiments in \cite{r437}
showed that the efficiency of the minimization greatly depended on the
deterministic placement of the interior points, with better results
obtained for these points placed sufficiently close to the boundary $S$ of
the obstacle $D$, but not very close to it. The current choice of a random
placement of the interior points $x_j$ reduced the variance in the
obtained results, and eliminated the need to provide a justified algorithm
 for their placement. The random choice of
these points distributes them in
 the entire interior of the obstacle,
rather than in a subset of it.

\subsection{Conclusions.}
For 3D obstacle  Rayleigh's hypothesis (conjecture) says that the
acoustic field
$u$ in the exterior of the obstacle $D$ is given by
the series convergent up to the boundary of $D$:
\be
u(x, \alpha) = e^{ik\alpha \cdot x} +
\sum^\infty_{\ell =0} A_\ell (\alpha)\psi_\ell, \quad
\psi_\ell:= Y_\ell (\alpha^\prime) h_\ell (kr), \quad
\alpha^\prime =
  \frac{x}{r}.  \end{equation}
While this conjecture (RC) is
false for many obstacles, it has been modified in \cite{r430} to obtain a
valid representation for the solution of
(\ref{s1_e1})-(\ref{s1_e2}). This representation (Theorem
\ref{s1_thm1}) is called the Modified Rayleigh Conjecture
(MRC), and is, in fact, not a conjecture, but a Theorem.

Can one use this approach to obtain solutions to various
scattering problems? A straightforward numerical
implementation of the MRC may fail, but, as we show here,
it can be efficiently implemented and allows one to
obtain accurate numerical solutions to  obstacle
scattering
problems.

The Random Multi-point
MRC algorithm was successfully applied to various 2D and 3D
obstacle scattering problems. This algorithm is a
significant improvement over previous MRC implementation
described in \cite{r437}.
The improvement is achieved by allowing the required
minimizations to be done iteratively, while the previous
methods were limited by the problem size constraints. In
\cite{r437}, such MRC method was presented, and it favorably
compared to the Boundary Integral Equation Method.

The Random Multi-point MRC has an additional attractive
feature, that it can easily treat obstacles with complicated
geometry (e.g.  edges and corners). Unlike the BIEM, it is
easily modified to treat different obstacle shapes.

Further research on  MRC algorithms is conducted. It is hoped that the MRC
 in its various implementation can emerge as a valuable and
efficient
alternative to more established methods.


\section{Support Function Method for inverse obstacle scattering problems.}

\subsection{Support Function Method (SFM)}

The Inverse Scattering Problem consists of finding the obstacle $D$ from
the Scattering Amplitude, or similarly observed data.
The Support Function Method (SFM) was originally developed in a
3-D setting in \cite{r56}, see also \cite{r190},  pp 94-99. It is used to approximately locate the
obstacle $D$. The method is derived using a high-frequency approximation
to the scattered field for smooth, strictly
 convex obstacles. It turns out that this inexpensive method also
provides a good localization of obstacles in the resonance
region of frequencies. If the obstacle is not convex, then the SFM
yields its convex hull.

One can restate the SFM  in a 2-D setting as follows (see \cite{r453}).
Let $D\subset\rc^2$ be a smooth and strictly convex obstacle with the
boundary $\Gamma$.
 Let $\nu(\by)$ be the
unique outward unit normal vector to $\Gamma$ at $\by\in\Gamma$.
Fix an incident direction $\a\in S^1$.
 Then the boundary $\Gamma$ can be decomposed into the following two parts:
\be
\Gamma_+=\{\by\in\Gamma\ :\ \nu(\by)\cdot
\a<0\}\,,\;\textnormal{and}\quad
\Gamma_-=\{\by\in\Gamma\ :\ \nu(\by)\cdot \a\geq 0\}\,,
\end{equation}
which are, correspondingly, the illuminated and the shadowed parts of the
boundary for the
chosen incident direction $\a$.

Given $\a\in S^1$, its {\bf specular point} ${\bf s_0}(\a)\in\Gamma_+$ is
defined from the condition:
\be
{\bf s_0}(\a)\cdot\a=\min_{{\bf s}\in\Gamma_+} {\bf s}\cdot\a
\end{equation}
Note that the equation of the tangent line to $\Gamma_+$ at ${\bf s_0}$
is
\be\label{s2_tand}
<x_1,x_2>\cdot\ \a={\bf s_0}(\a)\cdot \a\,,
\end{equation}
and
\be
\nu({\bf s_0}(\a))=-\a\,.
\end{equation}

The {\bf Support function} $d(\a)$ is defined by
\be
d(\a)={\bf s_0}(\a)\cdot \a\,.
\end{equation}

Thus $\abs{d(\a)}$ is the distance from the origin to the unique tangent
line to $\Gamma_+$ perpendicular to the incident vector $\a$.
Since the obstacle $D$ is assumed to be convex
\be\label{s2_suppd}
D=\cap_{\a\in S^1}\{\bx\in\rc^2\ :\ \bx\cdot\a\geq d(\a)\}\,.
\end{equation}
The boundary $\Gamma$ of $D$ is smooth, hence so is the Support Function. The
knowledge of this function allows one to reconstruct the boundary
$\Gamma$ using the following procedure.

Parametrize  unit vectors $\bl\in S^1$ by $\bl(t)=(\cos t,\sin t),\quad 0\leq t<2\pi$
and define
\be\label{s2_para2}
p(t)=d(\bl(t)),\quad 0\leq t<2\pi\,.
\end{equation}
Equation (\ref{s2_tand}) and the definition of the Support Function give
\be\label{s2_tan}
x_1\cos t+x_2\sin t=p(t)\,.
\end{equation}
Since $\Gamma$ is the envelope of its tangent lines, its equation can
be
found from (\ref{s2_tan}) and
\be
-x_1\sin t+x_2\cos t=p'(t)\,.
\end{equation}
Therefore
the parametric equations of the boundary $\Gamma$ are
\be\label{s2_pareq}
x_1(t)=p(t)\cos t-p'(t)\sin t,\quad
x_2(t)=p(t)\sin t+p'(t)\cos t\,.
\end{equation}

So, the question is how to construct the Support function $d(\bl),\;\bl\in S^1$ from the
knowledge of the Scattering Amplitude.
In 2-D the Scattering Amplitude is related to the total field
$u=u_0+v$ by

\be\label{s2_scat1}
A(\a', \a)=-\frac{e^{i\frac \pi 4}}{\sqrt{8\pi k}}
\int_\Gamma\frac{\partial u}{\partial\nu(\by)}\ e^{-ik\a'\cdot \by}\
ds(\by)\,.
\end{equation}

In the case of the "soft" boundary condition
(i.e. the pressure field satisfies the Dirichlet
boundary condition $u=0$) the Kirchhoff (high frequency) approximation
gives
\be
\frac{\partial u}{\partial \nu}=2\frac{\partial u_0}{\partial \nu}
\end{equation}
on the illuminated part $\Gamma_+$ of the boundary $\Gamma$, and
\be
\frac{\partial u}{\partial \nu}=0
\end{equation}
on the shadowed part $\Gamma_-$. Therefore, in this approximation,

\be\label{s2_scat2}
A(\a', \a)=-\frac{ike^{i\frac \pi 4}}{\sqrt{2\pi k}}
\int_{\Gamma_+}\a\cdot\nu(\by)\ e^{ik(\a-\a')\cdot \by}\
ds(\by)\,.
\end{equation}

Let $L$ be the length of $\Gamma_+$, and $\by=\by(\zeta),\; 0\leq\zeta\leq L$
be its arc length parametrization. Then
\be\label{s2_scat3}
A(\a', \a)=-\frac{i\sqrt{k}\ e^{i\frac \pi 4}}{\sqrt{2\pi}}
\int_0^L\a\cdot\nu(\by(\zeta))\ e^{ik(\a-\a')\cdot \by(\zeta)}
\ d\zeta\,.
\end{equation}

Let $\zeta_0\in[0,L]$ be such that ${\bf s}_0=\by(\zeta_0)$ is the specular point of
the unit vector $\bl$, where
\be\label{s2_scat31}
\bl=\frac{\a-\a'}{\abs{\a-\a'}}\,.
\end{equation}
Then
$\nu({\bf s}_0)=-\bl$, and $d(\bl)=\by(\zeta_0)\cdot\bl$.
Let
\[
\varphi(\zeta)=(\a-\a')\cdot \by(\zeta)\,.
\]
Then $\varphi(\zeta)=\bl\cdot\by(\zeta)\abs{\a-\a'}$.
Since $\nu({\bf s}_0)$ and $\by'(\zeta_0)$ are
orthogonal, one has
\[
\varphi'(\zeta_0)=\bl\cdot \by'(\zeta_0)\abs{\a-\a'}=0\,.
\]
Therefore, due to the strict convexity of $D$, $\zeta_0$ is also the
unique non--degenerate
stationary point of $\varphi(\zeta)$
 on the interval $[0,L]$, that is
$\varphi'(\zeta_0)=0$, and $\varphi''(\zeta_0)\not=0$.

According to the Stationary Phase method
\be\label{s2_scat4}
\int_0^L f(\zeta)e^{ik\varphi(\zeta)} d\zeta=f(\zeta_0)
\exp\left[ik\varphi(\zeta_0)+
\frac{i\pi}4 \frac{\varphi''(\zeta_0)}{\abs{\varphi''(\zeta_0)}}\right]
\sqrt{\frac
{2\pi}{k\abs{\varphi''(\zeta_0)}}}\left[1+O\left(\frac
1k\right)\right]\,,
\end{equation}
as $k\ra\infty$.

By the definition of the curvature $\kappa(\zeta_0)=\abs{\by''(\zeta_0)}$.
Therefore, from the collinearity of $\by''(\zeta_0)$ and $\bl$,
$\abs{\varphi''(\zeta_0)}=\abs{\a-\a'}\kappa(\zeta_0)$.
Finally, the strict convexity of $D$, and the definition of
$\varphi(\zeta)$, imply
that $\zeta_0$ is the unique point of minimum of $\varphi$ on
$[0,L]$, and
\be\label{s2_scat5}
\frac{\varphi''(\zeta_0)}{\abs{\varphi''(\zeta_0)}}=1\,.
\end{equation}
Using (\ref{s2_scat4})-(\ref{s2_scat5}),  expression
(\ref{s2_scat3}) becomes:

\be\label{s2_scat6}
A(\a',\a)=-\frac{\bl\cdot\a}{\sqrt{\abs{\a-\a'}\kappa(\zeta_0)}}
e^{ik(\a-\a')\cdot \by(\zeta_0)}
\left[1+O\left(\frac 1k\right)\right]\,, \quad
k\ra\infty\,.
\end{equation}

At the specular point one has
$\bl\cdot\a'=-\bl\cdot\a$. By the
definition $\a-\a'=\bl\abs{\a-\a'}$. Hence $\bl\cdot(\a-\a')=\abs{\a-
\a'}$ and $2\bl\cdot\a=\abs{\a-\a'}$. These equalities and $d(\bl)=
\by(\zeta_0)\cdot\bl\ $ give

\be\label{s2_scat7}
A(\a',\a)=-\frac 12\sqrt{\frac{\abs{\a-\a'}}{\kappa(\zeta_0)}}
e^{ik\abs{\a-\a'}d(\bl)}
\left[1+O\left(\frac 1k\right)\right]\,, \quad
k\ra\infty\,.
\end{equation}

Thus, the approximation
\be\label{s2_scat8}
A(\a',\a)\approx -\frac 12\sqrt{\frac{\abs{\a-\a'}}{\kappa(\zeta_0)}}
e^{ik\abs{\a-\a'}d(\bl)}
\end{equation}
can be used for an approximate recovery of the curvature and the support function
(modulo $2\pi/k\abs{\a-\a'}$)
 of the obstacle, provided one knows that the total field satisfies the
Dirichlet boundary condition.
 The uncertainty in the support function determination can be remedied
by using different combinations of vectors $\a$ and $\a'$ as described in the numerical results section.

Since it is also of interest to localize the obstacle in the case when the
boundary condition is not a priori known, one can modify the SFM as
shown in \cite{r463}, and obtain
\be\label{ss1_6}
A(\a',\a)\sim
\frac 12 \sqrt{\frac{\abs{\a-\a'}}{\kappa(\zeta_0)}}
\ e^{i(k\abs{\a-\a'}d(\bl)-2\gamma_0+\pi)},
\end{equation}
where
\[
\gamma_0=\arctan\frac kh,
\]
and
\[
\frac{\partial u}{\partial n}+hu=0
\]
along the boundary $\Gamma$ of the sought obstacle.

Now one can recover the Support Function $d(\bl)$ from (\ref{ss1_6}),
and the location of the obstacle.

\subsection{Numerical results for the Support Function Method.}

In the first numerical experiment the obstacle is the circle
\be
D=\{(x_1,x_2)\in\rc^2\ :\ (x_1-6)^2+(x_2-2)^2=1\}\,.
\end{equation}
It is reconstructed using the Support Function Method for two
frequencies in the resonance region: $k=1.0$, and $k=5.0$.
Table \ref{sfm_table3} shows how well the approximation (\ref{s2_scat8}) is satisfied
for various pairs of vectors $\a$ and $\a'$ all representing the same
vector $\bl=(1.0,0.0)$ according to (\ref{s2_scat31}). The Table shows
the ratios of the approximate Scattering Amplitude $A_a(\a',\a)$ defined
as the right hand side of the equation (\ref{s2_scat8}) to the exact
Scattering Amplitude $A(\a',\a)$. Note, that for a sphere of radius
$a$, centered at $\bx_0\in\rc^2$, one has
\be\label{s3_bfar}
A(\a',\a)=-\sqrt{\frac 2{\pi k}}\ e^{-i\frac{\pi}4}
e^{ik(\a-\a')\cdot \bx_0}
\sum_{l=-\infty}^\infty \frac
{J_l(ka)}{H^{(1)}_l(ka)}
\  e^{il(\theta-\beta)}\,,
\end{equation}
where $\a'=\bx/\abs{\bx}=e^{i\theta}$, and $\a=e^{i\beta}$.
Vectors $\a$ and $\a'$ are defined by their polar angles shown in
Table \ref{sfm_table3}.

\begin{table}[th]
\caption{Ratios of the approximate and the exact Scattering Amplitudes $A_a(\a',\a)/A(\a',\a)$
for $\bl=(1.0,0.0)$.\vspace*{1pt}}
{\footnotesize
\begin{tabular}{r r r r}
\hline
{} &{} &{} &{} \\[-1.5ex]
$\a'$ & $\a$ & $k=1.0\quad\ $ & $k=5.0\quad\ $ \\[1ex]
\hline
{} &{} &{} &{} \\[-1.5ex]
$\pi$ & $0$ &          0.88473   - 0.17487$i$ &0.98859   - 0.05846$i$\\[1ex]
$23\pi/24$ & $\pi/24$ &    0.88272   - 0.17696$i$ &0.98739   - 0.06006$i$\\[1ex]
$22\pi/24$ & $2\pi/24$ &   0.87602   - 0.18422$i$ & 0.98446   - 0.06459$i$\\[1ex]
$21\pi/24$ & $3\pi/24$ &   0.86182   - 0.19927$i$ &  0.97977   - 0.07432$i$\\[1ex]
$20\pi/24$ & $4\pi/24$ &   0.83290   - 0.22411$i$ &  0.96701   - 0.08873$i$\\[1ex]
$19\pi/24$ & $5\pi/24$ &   0.77723   - 0.25410$i$ &0.95311   - 0.10321$i$\\[1ex]
$18\pi/24$ & $6\pi/24$ &   0.68675   - 0.27130$i$ & 0.92330   - 0.14195$i$\\[1ex]
$17\pi/24$ & $7\pi/24$ &   0.57311   - 0.25360$i$ & 0.86457   - 0.14959$i$\\[1ex]
$16\pi/24$ & $8\pi/24$ &   0.46201   - 0.19894$i$ &0.81794   - 0.22900$i$\\[1ex]
$15\pi/24$ & $9\pi/24$ &   0.36677   - 0.12600$i$ & 0.61444   - 0.19014$i$\\[1ex]
$14\pi/24$ & $10\pi/24$ &  0.28169   - 0.05449$i$ & 0.57681   - 0.31075$i$\\[1ex]
$13\pi/24$ & $11\pi/24$ &  0.19019  + 0.00075$i$ & 0.14989   - 0.09479$i$\\[1ex]
$12\pi/24$ & $12\pi/24$ &  0.00000   + 0.00000$i$ & 0.00000   + 0.00000$i$\\[1ex]
\hline
\end{tabular}\label{sfm_table3} }
\end{table}

Table \ref{sfm_table3} shows that only vectors $\a$ close to the vector $\bl$ are
suitable for the Scattering Amplitude approximation. This shows the
practical importance of the backscattering data.
Any single combination of vectors $\a$ and $\a'$
representing $\bl$ is not sufficient to uniquely determine
the Support Function $d(\bl)$ from (\ref{s2_scat8}) because of the phase
uncertainty. However, one can remedy this by using more than one pair of vectors
$\a$ and $\a'$ as follows.

Let $\bl\in S^1$ be fixed. Let

\[
R(\bl)=\{\a\in S^1\ :\ \abs{\a\cdot\bl}>1/\sqrt{2}\}\,.
\]

Define $\Psi:\rc\ra \rc^+$ by
\[
\Psi(t)=\left\|\frac{A(\a',\a)}{\abs{A(\a',\a)}}+
e^{ik\abs{\a-\a'}t}\right\|^2_{L^2(R(\bl))}\,,
\]
where $\a'=\a'(\a)$ is defined by $\bl$ and $\a$ according to
(\ref{s2_scat31}), and the integration is done over
$\a\in R(\bl)$.

\begin{figure}[th]
\includegraphics*{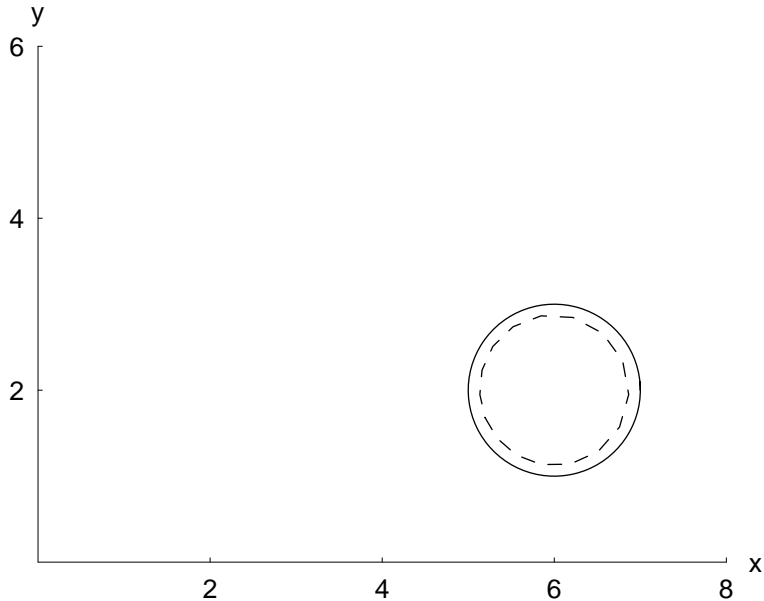}
\caption{Identified (dotted line), and the original (solid line) obstacle $D$ for $k=1.0$.}
\label{sfm_figure1}
\end{figure}

If the approximation (\ref{s2_scat8}) were exact for any $\a\in R(\bl)$,
then the value of
$\Psi(d(\bl))$ would be zero.
This justifies the use of the minimizer $t_0\in\rc$ of the function $\Psi(t)$ as an approximate value of the
Support Function $d(\bl)$. If the Support Function is known for
sufficiently many directions $\bl\in S^1$, the obstacle can be localized
using (\ref{s2_suppd}) or (\ref{s2_pareq}). The results of such a
localization for $k=1.0$ together with the original obstacle $D$ is shown
on Figure \ref{sfm_figure1}. For $k=5.0$ the identified obstacle is not shown, since it is practically the same as $D$.
 The only a priori assumption on $D$ was that it
was located inside the circle of radius $20$ with the center in the origin.
The Support Function was computed for 16 uniformly
distributed in $S^1$ vectors $\bl$.
The program run takes about 80 seconds on a 333 MHz PC.

In another numerical experiment we used $k=1.0$ and a kite-shaped obstacle. Its
boundary is described by
\be
{\bf r}(t)=(5.35+\cos t+0.65\cos 2t,\ 2.0+1.5\sin t),\quad 0\leq t<2\pi\,.
\end{equation}
Numerical experiments using the boundary integral equation method (BIEM)
for the direct scattering problem for
this obstacle centered in the origin are described in
\cite{coltonkress}, section 3.5. Again, the Dirichlet boundary
conditions were assumed.
We computed the scattering amplitude
for 120 directions $\a$
using the MRC method  with about 25\% performance
improvement over the BIEM, see
\cite{r437}.

\begin{figure}
\includegraphics*{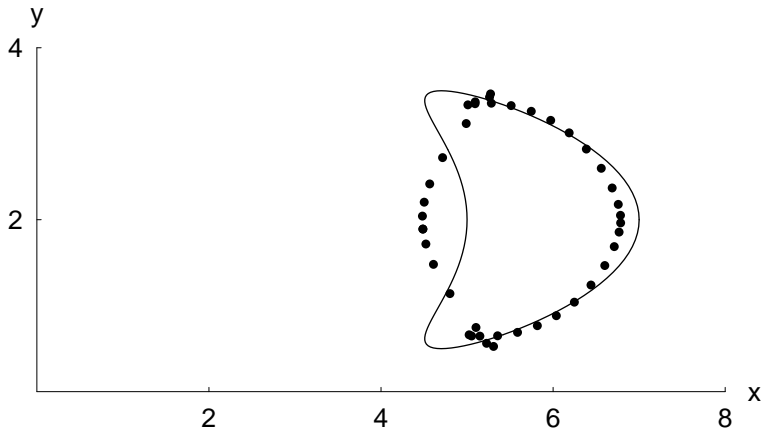}
\caption{Identified points and the original  obstacle $D$ (solid line);
 $k=1.0$.}\label{mrc_figure3}
\end{figure}

The Support Function Method (SFM) was used to identify the obstacle $D$ from
the synthetic scattering amplitude with no noise added.
 The only a priori assumption on $D$ was that it
was located inside the circle of radius $20$ with the center in the origin.
The Support Function was computed for 40 uniformly
distributed in $S^1$ vectors $\bl$ in about 10 seconds on a 333 MHz PC. The results
of the identification are shown in Figure \ref{mrc_figure3}. The original obstacle is
the solid line. The points were identified according to (\ref{s2_pareq}).
As expected, the method recovers
the convex part of the boundary $\Gamma$, and fails for the concave part.
The same experiment but with $k=5.0$ achieves a perfect identification of the convex
part of the boundary. In each case the convex part of the obstacle was successfully localized.
Further improvements in the obstacle localization using the MRC method are suggested
in \cite{r430}, and in the next section.

For the identification of obstacles with unknown boundary conditions let
\[
A(t)=A(\a',\a)=\abs{A(t)}e^{i\psi(t)}
\]
where, given $t$, the vectors $\a$ and $\a'$ are chosen as above, and the phase function
$\psi(t),\;\sqrt{2}<t\leq 2$ is continuous.
Similarly, let
$A_a(t),\;\psi_a(t)$ be the approximate scattering amplitude
 and its phase defined by formula (\ref{ss1_6}).

If the approximation (\ref{ss1_6}) were exact for any $\a\in R(\bl)$,
then the value of
\[
\abs{\psi_a(t)-ktd(\bl)+2\gamma_0-\pi}
\]
 would be a multiple of $2\pi$.

This justifies the following algorithm for the determination of the Support Function
$d(\bl)$:

Use a linear regression to find the approximation
\[
\psi(t)\approx C_1t+C_2
\]
 on the interval $\sqrt{2}<t\leq 2$. Then
\be\label{s1_appra}
d(\bl)=\frac {C_1} k\,.
\end{equation}
Also
\[
h=-k\tan\frac {C_2} 2\,.
\]
However, the formula for $h$ did not work well numerically. It could only
determine if the boundary conditions were or were not of the Dirichlet type.
Table \ref{sfm_table4} shows that
the algorithm based on (\ref{s1_appra})
was successful in the identification of the circle of radius
$1.0$ centered in the origin for various values of $h$ with no a priori
assumptions on the boundary conditions.
For this circle the Support Function $d(\bl)=-1.0$ for any direction $\bl$.

\begin{table}
\caption{Identified values of the Support Function for the circle of radius $1.0$ at $k=3.0$.}

\begin{tabular}{ r  r  r }

\hline
$h$ & Identified $d(\bl)$ & Actual $d(\bl)$ \\

\hline
  0.01   &  -0.9006 &      -1.00 \\
  0.10   &   -0.9191 &     -1.00 \\
  0.50   &  -1.0072  &    -1.00 \\
  1.00   &  -1.0730  &    -1.00 \\
  2.00   &  -0.9305  &    -1.00 \\
 5.00    &   -1.3479  &    -1.00 \\
 10.00   &  -1.1693   &   -1.00 \\
 100.00  &  -1.0801  &     -1.00 \\
\hline
\end{tabular}\label{sfm_table4}

\end{table}

\section{Analysis of a Linear Sampling method.}

During the last decade many papers were published, in which the obstacle identification
methods were based on a numerical verification of the inclusion of some function
$f:=f(\alpha,z), \, z\in \R^3, \, \alpha \in S^2,$ in
the range $R(B)$ of a certain operator $B$. Examples of such methods
include  \cite{ck}, \cite{cc},\cite{Kir98}. However, one can show that
 the methods proposed in the
above papers have essential difficulties, see \cite{r462}. Although it is true that
$f\not\in R(B)$ when $z\not\in D$, it turns out that in any
neighborhood of $f$ there are elements from $R(B)$. Also, although
$f\in R(B)$ when $z\in D$, there are elements in every neighborhood of
$f$ which do not belong to $R(B)$ even if $z\in D$.
Therefore it is quite difficult to construct a
stable numerical method for the identification of $D$ based on the verification of the
inclusions $f\not\in R(B)$, and $f\in R(B)$. Some published numerical
 results were intended to show that the method based on the above idea works
practically, but it is not clear how these conclusions were obtained.

Let us introduce some {\it notations} : $N(B)$ and $R(B)$
are, respectively, the null-space and the range of a linear operator
$B$, $D\in \R^3$ is a
bounded
domain (obstacle) with a smooth boundary $S$, $D'=\R^3 \setminus D,$
$u_0=e^{ik\alpha \cdot x}$, $k=const>0$, $\alpha\in S^2$ is a unit vector,
$N$ is the unit normal to $S$ pointing into $D'$, $g=g(x,y,k):=g(|x-y|):=
\frac {e^{ik|x-y|}}{4 \pi |x-y|}$, $f:=e^{-ik \alpha' \cdot z}$, where
$z\in
\R^3$ and $\a' \in S^2$, $\a':=xr^{-1}$,
$r=|x|$, $u=u(x,\alpha, k)$ is the scattering solution:
\be\label{lsm_1}
(\Delta +k^2)u=0 \quad in \quad D', u|_{S}=0,
\end{equation}
\be\label{lsm_2}
u=u_0+v,\quad v=A(\a',\alpha,k)e^{ikr}r^{-1} +o(r^{-1}), \quad as \quad
r\to \infty,
\end{equation}
where $A:=A(\a',\alpha,k)$ is called the scattering amplitude,
corresponding
to the obstacle $D$ and the Dirichlet boundary condition.
Let $G=G(x,y,k)$ be the resolvent kernel of the Dirichlet Laplacian in
$D'$:
\be\label{lsm_3}
(\Delta +k^2)G=-\d(x-y) \quad in \quad D', G|_{S}=0,
\end{equation}
and $G$ satisfies the outgoing radiation condition.

If
\be\label{lsm_4}
(\Delta +k^2)w=0 \quad in \quad D', w|_{S}=h,
\end{equation}
and $w$ satisfies the radiation condition, then (\cite{r190})
one has
\be\label{lsm_5}
w(x)=\int_S G_N(x,s)h(s)ds,\quad w=A(\a',k)e^{ikr}r^{-1} +o(r^{-1}),
\end{equation}
as $r\to\infty$, and $xr^{-1}=\a'$.
We write $A(\a')$ for $A(\a',k)$, and
\be\label{lsm_6}
A(\a'):=Bh:= \frac 1 {4\pi}\int_S u_N(s,-\a')h(s)ds,
\end{equation}
as follows from Ramm's lemma:

{\bf Lemma 1.\,\,} (\cite{r190}, p.46) {\it One has:
\be\label{lsm_7}
G(x,y,k)= g(r)u(y,-\a',k) +o(r^{-1}),
\quad as \quad r=|x|\to
\infty, \quad xr^{-1}=\a',
\end{equation}
where $u$ is the scattering solution of (\ref{lsm_1})-(\ref{lsm_2})}.

 One can write the scattering amplitude  as:
\be\label{lsm_8}
A(\a',\alpha,k)=-\frac 1 {4\pi}\int_S u_N(s,-\a')e^{ik\alpha \cdot s}ds.
\end{equation}
The following claim follows easily from the results in \cite{r190}, \cite{rammb2} (cf
\cite{Kir98}):

{\it Claim:}  $f:=e^{-ik \alpha' \cdot z}\in R(B)$ if and only if
$z\in D$.

{\it Proof:} If $e^{-ik \alpha' \cdot z}=Bh$, then Lemma 1
 and (12.6) imply
\[g(y,z)=\int_SG_N(s,y)hds\quad \text{for}\quad |y|>|z|\,.
\]
Thus $z\in D$, because otherwise one gets a contradiction: $\lim_{y\to z} g(y,z)=\infty$
if $z\in \overline {D'}$ , while $\lim_{y\to z}\int_SG_N(s,y)hds<\infty $
if $z\in \overline {D'}$. Conversely, if $z\in D$, then Green's formula yields
$g(y,z)=\int_SG_N(s,y)g(s,z)ds$. Taking $|y|\to \infty, \, \frac y {|y|}=\alpha'$, and using
Lemma 1, one gets  $e^{-ik \alpha' \cdot z}=Bh,$ where $h=g(s,z)$. The claim is proved. $\Box$

Consider $B: L^2(S)\to L^2(S^2)$, and $A: L^2(S^2)\to L^2(S^2)$,
where $B$ is defined in (\ref{lsm_6}) and $Aq:=\int_{S^2}A(\a',\alpha)q(\alpha)
d\alpha.$ Then one proves (see \cite{r462}):

{\bf Theorem 1.} {\it The ranges $R(B)$ and $R(A)$ are dense in $L^2(S^2)$
}

{\bf Remark 1.} In \cite{ck} the 2D inverse obstacle scattering
problem is considered. It is proposed to solve the equation
(1.9) in \cite{ck}:
\be\label{lsm_9}
\int_{S^1}A(\alpha,\b)\g d\b=e^{-ik\alpha \cdot z},
\end{equation}
where $A$ is the scattering amplitude at a fixed $k>0$, $S^1$ is the unit
circle, $\alpha\in S^1$, and $z$ is a point on $\R^2$. If $\g=\g(\b, z)$ is
found, the
boundary
$S$ of the obstacle is to be found by finding those $z$ for which
$||\g||:=||\g(\b,z)||_{L^2(S^1)}$ is maximal.
Assuming that $k^2$ is not a Dirichlet or Neumann eigenvalue
of the Laplacian in $D$, that $D$ is a smooth, bounded, simply
connected  domain, the authors state Theorem 2.1 \cite{ck}, p.386,
which says that for every $\epsilon>0$ there exists a function $\g\in
L^2(S^1),$ such that
\be\label{lsm_10}
\lim_{z\to S}||\g(\b,z)||=\infty,
\end{equation}
and (see \cite{ck}, p.386),
\be\label{lsm_11}
||\int_{S^1}A(\alpha,\b)\g d\b-e^{-ik\alpha \cdot z}||<\epsilon.
\end{equation}

 There are several questions concerning the proposed method.

First, equation (\ref{lsm_9}), in general, is not solvable. The authors propose to
solve it approximately, by a regularization method. The regularization
method applies for stable solution of solvable
ill-posed equations (with exact or noisy data). If equation (\ref{lsm_9})
is not solvable, it is not clear what numerical "solution"
one seeks by a regularization method.

Secondly, since the kernel of the integral operator in (\ref{lsm_9}) is
smooth, one can always find, for any $z\in \R^2$, infinitely many $\g$ with
arbitrary large $||\g||$, such that (\ref{lsm_11}) holds.
 Therefore it is not clear how and why, using (\ref{lsm_10}), one can find $S$
numerically by the proposed method.

A numerical implementation of the Linear Sampling Method (LSM) suggested in \cite{ck}
consists of solving a discretized version of (\ref{lsm_9})
\be\label{lsm_12}
F\bg=\bff,
\end{equation}
where $F=\{A{\a_i,\b_j}\},\; i=1,...,N$, $j=1,...,N$ be a square matrix
formed by the measurements of the scattering amplitude for $N$ incoming,
and $N$ outgoing directions. In 2-D the vector $\bff$ is formed by
$$
\bff_n=\frac{e^{i\frac \pi 4}}{\sqrt{8\pi k}}e^{-ik\a_n\cdot z},\quad
n=1,...,N,
$$
see \cite{brandfass} for details.

Denote the Singular Value Decomposition of the far field operator by
$F=USV^H$. Let $s_n$ be the singular values of $F$, $\rho=U^H\bff$, and $\mu=V^H \bff$. Then the norm of
the sought function $g$ is given by
\be\label{lsm_13}
\|\g\|^2=\sum^N_{n=1}\frac{|\rho_n|^2}{s_n^2}.
\end{equation}
A different LSM is suggested by A. Kirsch in \cite{Kir98}. In it one solves
\be\label{lsm_14}
(F^*F)^{1/4}\bg=\bff
\end{equation}
instead of (\ref{lsm_12}). The corresponding expression for the norm of
$\g$ is
\be\label{lsm_15}
\|\g\|^2=\sum^N_{n=1}\frac{|\mu_n|^2}{s_n}.
\end{equation}
A detailed numerical comparison of the two LSMs and the linearized
tomographic inverse scattering is given in \cite{brandfass}.

\begin{figure}
\includegraphics*{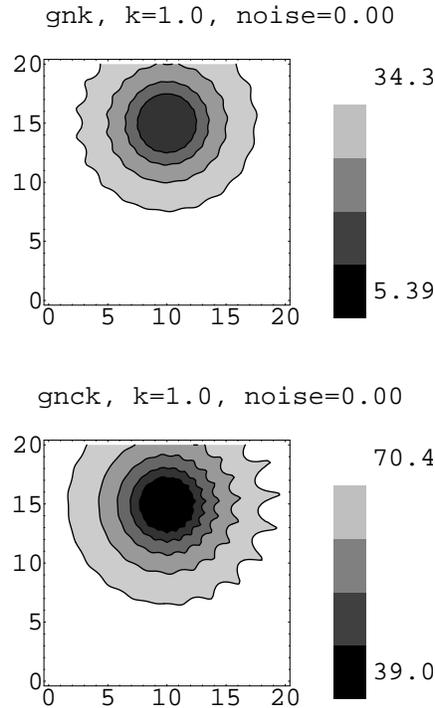}
\caption{Identification of a circle at $k=1.0$.}\label{lsm_figure4}
\end{figure}

\begin{figure}
\includegraphics*{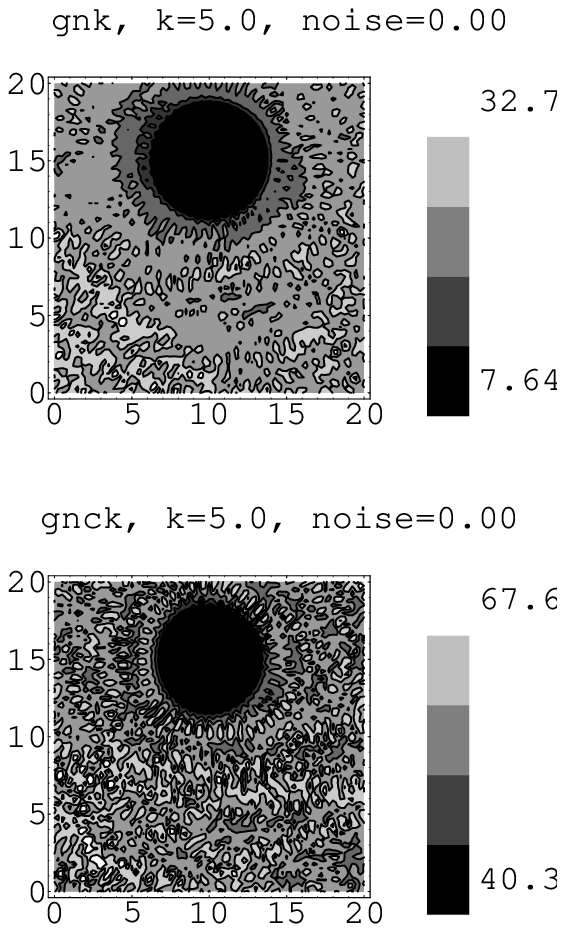}
\caption{Identification of a circle at $k=5.0$.}\label{lsm_figure5}
\end{figure}

 The conclusions of \cite{brandfass}, as well as of our own numerical
experiments are that the method of Kirsch (\ref{lsm_14}) gives a better,
but a comparable identification, than (\ref{lsm_12}). The identification
is significantly deteriorating if the scattering amplitude is available
only for a limited aperture, or the data are corrupted by noise. Also, the points with the
{\it smallest} values of the $\|\g\|$ are the best in locating the
inclusion, and not the {\it largest} one, as required by the theory in \cite{Kir98}
and in \cite{ck}.
 In Figures \ref{lsm_figure4} and \ref{lsm_figure5} the
implementation of the
Colton-Kirsch LSM (\ref{lsm_13}) is denoted by $gnck$, and of the
Kirsch method (\ref{lsm_15}) by $gnk$. The Figures show a contour plot of the logarithm
of the $\|\g\|$. In all the  cases the original
obstacle was the circle
of radius $1.0$ centered at the point $(10.0,\ 15.0)$. A similar circular obstacle that
was identified by the Support Function Method (SFM) is discussed in Section 10.
Note that the actual radius of the circle is $1.0$, but it cannot be seen from the
LSM identification. The LSM does not require any knowledge of the
boundary conditions on the obstacle. The use of the SFM for unknown
boundary conditions is discussed in the previous section. The LSM identification
was performed for the scattering amplitude of the circle computed
analytically with no noise added. In all the experiments the value for the parameter
 $N$ was chosen to be 128.

\end{document}